\documentclass[fleqn,11pt]{article}
\usepackage{amssymb,amsmath,amsfonts}
\usepackage{color}
\usepackage{graphicx}
\usepackage{latexsym}
\usepackage{amsmath}
\usepackage{amssymb}
\usepackage[dvips]{epsfig}
\usepackage[margin=1in]{geometry}

\numberwithin{equation}{section}
\newtheorem{thm}{Theorem}[section]
\newtheorem{coro}[thm]{Corollary}
\newtheorem{lem}[thm]{Lemma}

\newtheorem{rmk}[thm]{Remark}

\newcommand{\da}{\delta}

\newcommand{\ga}{\gamma}

\newcommand{\iy}{\infty}

\newcommand{\q}{\quad}
\newcommand{\e}{\epsilon}
\newcommand{\lt}{\left}
\newcommand{\rt}{\right}
\newcommand{\be}{\begin{equation}}
\newcommand{\bs}{\begin{split}}
\newcommand{\es}{\end{split}}
\newcommand{\ee}{\end{equation}}
\newcommand{\bee}{\begin{equation*}}
\newcommand{\eee}{\end{equation*}}

\newcommand{\ef}{\eqref}

\begin{document}
\begin{center}
\large{ \bf Low Mach and low Froude number limit for vacuum free boundary problem of all-time classical solutions of 1-D compressible Navier-Stokes equations}
\end{center}
\begin{center}
{Yaobin Ou\footnote{Email: ou@ruc.edu.cn, ou.yaobin@gmail.com}}
\end{center}

\begin{center}
{\small School of Mathematics, Renmin University of China,\\ Beijing 100872, P.R. China}
\end{center}

\abstract  In this paper, we study the low Mach and Froude number limit for the all-time classical solution of a fluid-vacuum free boundary problem  of  one-dimensional compressible Navier-Stokes equations. No smallness of initial data for the existence of all-time solutions are supposed.  The uniform estimates of solutions with respect to the Mach number and the Froude number are established for all the time, in particular for high order derivatives of the pressure, which is a novelty in contrast to previous results. The cases of ``ill-prepared" initial data and ``well-prepared" initial data are both discussed. It is interesting to see, either both the Mach number and the Froude number vanish, or the time goes to infinity, the limiting functions are the same, that is, the steady state. The main difficulty is that, the system is degenerate near the free boundary and contains singular terms. This result can be viewed as the first one on the low Mach and Froude numbers limit for free boundary problems. At the same time, we also establish the all-time existence of the classical solution with sharp convergent rates to the steady state, while previous results are only concerned with the weak or strong solutions.


\section{Introduction}
The motion of isentropic viscous compressible fluids in a gravity field is governed by the following compressible Navier-Stokes equations:
\begin{equation} \label{NS1}\rho_t+\textrm{div}(\rho u)=0,
\end{equation}
\begin{equation}\label{NS2}
(\rho u)_t+\textrm{div}(\rho u\otimes
u)-\textrm{div}\mathbb{S}(u)+\nabla p=-\rho g,
\end{equation}
where $\rho,u=(u^1,...,u^n),p=\rho^\gamma$ denote density, velocity and
pressure,  respectively. The constants $\mu,\lambda$
are viscous coefficients with $\mu>0$, $\lambda+\frac{2}{n}\mu\ge0$,
$
\mathbb{S}(u)=2\mu D(u)+\lambda{\rm div} \textrm{u}I
$
is the viscous stress tensor, with $D(u)=\frac{1}{2}(\nabla u+\nabla
u^t)$, and $g>0$ is the gravitational constant.   For the air at standard conditions,
$\gamma$ is approximately $1.4$.

 We recover the  equations with non-dimensional numbers by dividing each physical variable by
its reference state (for instance, the mean value) as follows:
$$\rho=\frac{\tilde{\rho}}{r},\;p=\frac{\tilde{p}}{\bar{p}},
 \;u=\frac{\tilde{u}}{\bar{u}},\;g=\frac{\tilde{g}}{\bar{g}}
$$
$$x=\frac{\tilde{x}}{L},\;t=\frac{\tilde{t}}{L/\bar{u}},
\;\mu=\frac{\tilde{\mu}}{\nu},\;\lambda=\frac{\tilde{\lambda}}{\nu},
$$
where $r,\bar{p},\bar{u},L,\nu$ are the reference
density, pressure, velocity,  length, and viscosity
coefficient, respectively. After dropping all the tildes in the
system and taking all other dimensionless quantities to be 1, we
obtain the compressible Navier-Stokes equations in dimensionless form as  follows (see also \cite{FMNS,Masmoudi}):
\begin{equation} \label{CNS1}\rho_t+\textrm{div}(\rho  u )=0,
\end{equation}
\begin{equation}\label{CNS2}
(\rho  u )_t+\textrm{div}(\rho  u \otimes
u )-\textrm{div}\mathbb{S}(u )+\frac{1}{\rm {Ma}^2}\nabla p =-\frac{1}{\rm {Fr}^2}\rho  g ,
\end{equation}
where the non-dimensional constants \be {\rm Ma}=\frac{\bar{u}}{\sqrt{\gamma\bar{p}/\bar{\rho}}},\quad {\rm Fr}=\frac{\bar u}{\bar g L}
\ee are the (reference) Mach number and Froude number, respectively.

Suppose ${\rm Ma=Fr=\epsilon}$ for simplicity, and
$(\rho^\epsilon,u^\epsilon)$ is the solution to
\eqref{CNS1}-\eqref{CNS2} for any $\epsilon>0$. Then as $\epsilon\to 0$, the sequence
\begin{equation*} \rho^\epsilon \to \bar \rho,\;u^\epsilon\to w,\;\epsilon^{-2}(\nabla
p^\epsilon+\rho^\epsilon g) \to \bar\rho\nabla\pi,
\end{equation*} where $(w,\pi)$ is a solution to the following anelastic system in  $\mathbb{R}^n$ or $\mathbb{T}^n$(see \cite{FMNS,Masmoudi}):
\begin{equation}\label{limit1}
\textrm{div}(\bar\rho w)=0,
\end{equation}
\begin{equation}\label{limit2}
(\bar\rho w)_t+\textrm{div}(\rho w\otimes w)-\textrm{div}\mathbb{S}(w)+\bar\rho\nabla \pi=0,
\end{equation}
and $(\bar\rho,0)$ is a solution to the steady problem of \eqref{limit1}-\eqref{limit2}.

The study of low Mach number limit which only vanishes the Mach number ${\rm Ma}$ (with 1/Fr=1 or 0) began with the pioneering work by Klainerman and Majda \cite{KM2} for the local strong solution with ``well-prepared'' initial data in the sense that the divergence of initial velocity is small in high order norms. Significant progresses have been achieved for global solutions in an arbitrary given time interval with ``ill-prepared'' initial data, which started with the work of  Lions and Masmoudi \cite{LM}, and continued by Danchin \cite{DAN1,DAN2}, Desjardins et al. \cite{DG,DGLM},  and the references cited therein. In these papers, the techniques for analyzing the acoustic waves, developed by Schochet \cite{SC1} and Grenier \cite{Grenier}, were applied and extended. The all-time existence and incompressible limit of  solutions to isentropic Navier-Stokes equations are studied by  Hagstrom and Lorenz \cite{HL}, Hoff \cite{Hoff}, Ou \cite{OU1}, for instance. In the presence of heat conduction or other physical factors (e.g., magnetic fields), rich phenomena occur in the process of the singular limit. One may refer to \cite{AL2,BDGL,DJO,FN,FN2,HW2,JJLX,JO,JLL} and the references therein for the contributions on related problems. For the free boundary problems, interesting results  on the formal or rigorous verification on the low Mach number limit have been obtained in \cite{Ali,FKNNS}, respectively.

In the regime that both the Mach number and the Froude number are small, the flow is strongly stratified. In particular, in the modeling of atmospheric flows (for instance, Ogura and Philips \cite{OP}), the high frequency sound waves are considered to be of low importance. However, it is difficult to verify the limit as both the Mach number and the Froude number tend to zero, since the process is singular. The rigorous mathematical studies of low Mach and Froude number limit were proposed by Masmoudi \cite{Masmoudi} for the case of bounded domains, and Feireisl et al. \cite{FMNS} for the case of periodic domains independently, and later extended to the situations of  non-isentropic flows with purely transported entropy (e.g., \cite{FKNZ}). However, the investigation of this singular limit for free boundary problems is still open, due to the lack of the corresponding well-posedness theory.

In this paper, we study the low Mach and Froude number limit of  the vacuum free  boundary problem  of one-dimensional   compressible Navier-Stokes equations \ef{CNS1}-\ef{CNS2} with large initial data, as well as the all-time existence of the classical solutions when the fluid connects to the vacuum continuously. Note that in case of vacuum free boundary problem, the all-time existence of classical solution with large initial data and decay rates is not available till now, even though in the one-dimension case. The main reason is due to the degeneracy near the free boundary. For the sake of simplicity and well understanding on the behaviors of solutions, we consider the one-dimensional case which reads as follows:
\be\label{1}\lt\{\begin{split}
&\rho_t^\e + \lt(\rho^\e u^\e\rt)_\eta =0         & {\rm in} &\ \  \cup_{\{t\ge 0\}}(0, {\Gamma^\e}(t)),\\
&(\rho^\e u^\e)_t +(\rho^\e (u^\e)^2)_\eta  +\frac{1}{\epsilon^2}p_\eta^\e = \mu u_{\eta\eta}^\e-\frac{1}{\epsilon^2}\rho^\e g   & {\rm in } &\ \ \cup_{\{t\ge0\}}(0, {\Gamma^\e}(t)),\\
& \rho^\e>0                                  & {\rm in } &\ \ \cup_{\{t\ge 0\}}[0, {\Gamma^\e}(t)),\\
& u^\e(0,t)=0, \  \    (\frac{1}{\epsilon^2}p^\e-\mu u_\eta^\e)(\Gamma^\e(t),t)=0 ,\\
      &  \dot{\Gamma^\e}(t) = u^\e  (\Gamma^\e(t),t),   \quad         t>0, \\
&(\rho^\e, u^\e)(x,0)=(\rho_0^\e, u_0^\e)                 & {\rm on} & \ \ (0,l_0).
 \end{split}\rt.\ee
Here $(\eta,t)\in \mathbb{R}\times [0,\iy)$,  $\rho^\e$, $u^\e $ and $p^\e$ denote the space and time variable, density, velocity and pressure, respectively;  $\Gamma^\e(t)$ and $ \dot{\Gamma^\e}(t) $ signify, respectively, the moving boundary between the fluid and the vacuum, and  velocity of $\Gamma^\e(t)$;  $g$ and $l_0$ are  positive constants;   $\rho^\e g$ appearing on $\ef{1}_{2}$  is the gravitational force.
The equation of state is given by
\be\label{2} p^\e=(\rho^\e)^\gamma,\;\gamma>1   \ee
where $\gamma$ is the adiabatic component.   The constants $\mu>0$ are set to be unity for convenience. In \ef{1}, the Mach number and the Froude number are both proportional to $\epsilon\in(0,1],$ for simplicity.

The initial density is assumed to satisfy the following conditions (see also \cite{LXZ}):
\be\label{initialdata}
\rho_0^\e>0 \ \ {\rm on} \ \ [0, l_0), \ \  \rho_0^\e(l_0)=0, \ \  -\infty< \lt((\rho_{0 }^\e)^{\ga-1}\rt)_x(l_0) <0 \ \ {\rm and} \ \  \int_0^{l_0}\rho_0^\e(x)dx=M.
\ee
It follows that $(\rho_0^\e)^{\ga-1} =O(1) (  l_0- x)$ as $x$ close to $l_0$.  From the conservation of mass, i.e. $\ef{1}_1$, we have
$$\int_0^{\Gamma^\e(t)} \rho^\e(x,t) dx = \int_0^{l_0}\rho_0^\e(x)dx=M. $$
The steady solution $(\rho, u)=(\bar\rho, 0)$ of \ef{1} with the same mass is defined uniquely by
\be\label{station}
 (\bar\rho^\gamma)_x=- \bar\rho g \ \ {\rm and} \ \   \int_0^\iy \bar\rho (x)dx=M,
\ee
which yields
\be\label{oo}
\bar\rho^{\gamma-1}(x)=  \frac{g(\gamma-1)}{\gamma}(\bar l-x),  \ \ {\rm for} \ \  0\le x\le \bar l:= \frac{\ga}{(\ga-1)g} (Mg)^{\frac{\ga-1}{\ga}},
\ee
and for $x> \bar l$, we set $\bar\rho(x)=0$.

We remark that in the rigorous derivation of the low Mach and Froude numbers limit for the one-dimensional free boundary problem, the process  and the limiting equations are  somewhat different from the cases of $\mathbb{T}^n$ or $\mathbb{R}^n$ or bounded domains   since the problem is treated in the Lagrangian trajectory coordinate, but not the Eulerian coordinate. The convergence is established in the sense of particle path (see \ef{conv2}), since the domain of fluid always changes. However, we find that the limiting solution is also $(\bar\rho(x),0)$, which is exactly a special solution to \ef{limit1}-\ef{limit2}.

The free boundary problems concerning the compressible fluids and the vacuum,   arise in many important
physical phenomena, such as gaseous stars and planets, shallow water waves, and so on. There have been some significant results on local results for isentropic inviscid flows ($\mu=0$, cf. \cite{coutlandshkoller,jm}) and viscous flows (cf. \cite{LiXY,YZ}). For the global-in-time theory, important progresses are achieved on weak solutions for the compressible viscous models with either the constant viscosity or density-dependent viscosities.  In case that the fluid connects to vacuum continuously, Okada (\cite{Okada}) showed the global existence of weak solutions to \ef{1}  and the large time asymptotic behaviors of solutions (without convergent rates). Fang and Zhang (\cite{fangzhang}), Zhu and Zi (\cite{zhuzi}), Ou and Zeng (\cite{OZ}) investigated the existence and large time behaviors of the global weak/strong solutions to the one-dimensional Navier-Stokes equations  with density-dependent viscosity ($\mu\approx \rho^\theta$)  and small data, while the large time stability of radially symmetric  strong solutions to the viscous gaseous stars problem was studied by Luo, Xin and Zeng (\cite{LXZ}), with the smallness assumption on the initial data. Recently, the global existence of large classical solutions for  compressible Navier-Stokes equations are established in \cite{zeng,ou}, without discussing the large-time behaviors of the free boundary and the solution.

In the theory of free boundary problems  for nonlinear partial differential equations, it is very important to study not only the global existence of  smooth solutions  with large initial data, but also  various behaviors of solutions.  In this paper, we  study the low Mach and Froude numbers limit and large-time asymptotic stability of the all-time classical solution to the vacuum free boundary problem \ef{1}. The uniform estimates with respect to the physical parameter $\epsilon\in (0,1]$ and the time $t\in[0,+\infty)$ are established; moreover, as $\e\to 0$ or $t\to\infty$, the solutions will converge to the same steady state solution with sharp rates. To the best of our knowledge, this is the first result on the low Mach and low Froude number limit for free boundary problems.

 Previous results on the low Mach and Froude numbers limit (see \cite{FKNZ,FMNS,FNP,Masmoudi}, for instance) mainly focus on the weak solutions of initial value problems or initial-boundary value problems in multi-dimensions.
 However, the singular limit of free boundary problems was not yet studied.
 In this paper, we are interested in the one for  the one-dimensional vacuum free boundary problems. In particular, the uniform estimates in $\epsilon\in (0,1]$ and $t\in[0,+\infty)$ for certain high order derivatives and the large-time stability are established  in the Lagrangian coordinate, which are the ingredients of this work. These properties are mostly due to the fine structure of the system in one spatial dimension where the weighted energy estimates can be derived inductively by using the Hardy-type inequalities. The methods also work for other kinds of singular limits of 1-D free boundary problems, however, they cannot be applied to the multi-dimensional free boundary problems directly, since  the geometry of the free boundary should be handled. This singular limit is also different from the incompressible limit in \cite{FKNNS} since the limiting density here is inhomogeneous, and the framework is for 1-D all-time classical solutions instead of multi-dimensional global weak solutions. Moreover, in contrast to previous stability results for constant viscosity and fixed $\epsilon$ (cf. \cite{Okada}), the high regularity of solutions and sharp convergent rates are obtained in this paper. Our methods can also be extended to handling the cases of density-dependent viscosity $\mu=\rho^\theta,\;\theta>0$ (cf. \cite{fangzhang,OZ,zhuzi}) in the sense that the additional high order  stability and the uniform estimates in $\e$ can be achieved. Note that the smallness assumption on the initial data is required in the situation of $\mu=\rho^\theta$ due to the degeneracy of the diffusion term  and the different boundary condition, which is a restriction of our technique.

The main difficulties in this result lie  in the degeneracy and singualrity of the system, since the density vanishes at the free boundary and large terms of the order $1/\e^2$ appear in the momentum equation. Thus standard methods for parabolic equations and the frame work for anti-symmetric singular differential operators (cf. \cite{KM2}) do not apply to get the uniform estimates with respect to $\epsilon\in (0,1]$ and $t\in[0,+\infty)$. Motivated by \cite{LXZ,OZ}, we overcome these troubles by converting the system into a degenerate parabolic equation with good structure by the Lagrangian trajectory. The key for the uniform estimates is to show the uniform upper and lower bounds in $\epsilon\in (0,1]$ and $t\in[0,+\infty)$ for the deformation variable $\eta_x$, which is indeed the Jacobian between the original coordinate and the resulting Lagrangian coordinate. We also introduce a new sort of singular  multipliers to derive the weighted energy estimates up to third order spatial derivatives, and then achieve the non-weighted estimates by a refined version of the Hardy inequality. Moreover, we  can derive the sharp decay rates of  the energy, by carefully choosing suitable weights and multipliers, applying the interpolation between lower-order energies and higher-order energies and  using delicately  Hardy's inequality.
The advantage  of our approach is that we can prove the uniform smoothness
and large time convergence of solutions with the detailed convergence rates simultaneously, without any smallness assumption.

 Finally, we remark that when the density is always bounded below by a positive constant, many important results are achieved during the past twenty years for the local or global solutions in isentropic regime or non-isentropic regime, either for the cases with constant viscosity or the
cases with density-dependent viscosity.  Please refer to \cite{GLX,JXZ,LiXY, QY, Wang,zhu1} and the references therein for details.

The rest of this paper is organized as follows. In Section 2, we  reformulate the free boundary problem \ef{1} by the Lagrangian particle path formulation and state the main results of this paper. In Section 3, we give some preliminaries of this paper. Next, in Section 4, we establish the uniform estimates and global existence for the solutions with ``ill-prepared'' initial data in $\epsilon\in (0,1]$ and $t\in[0,+\infty)$; we also obtain sharp convergent rates to the steady state solution as $t\to\infty$. In Section 5, we show the singular limit of the global solutions with ``ill-prepared'' initial data as both the Mach number and the Froude number  tend to zero. Finally, in Section 6, we show the local existence and uniqueness of  classical solutions to the problem \ef{1}.

\vskip 3mm
\noindent{\bf Notations:}

1) Throughout the rest of paper,   $C$  will denote a general positive constant which  only depends on the parameters of the problem, e.g.
$\ga$, but  not  on the initial data.  Also we use $C(\beta)$ to denote  a  general positive constant depending on the quantity $\beta$.

2) We will employ the notation $a\lesssim b$ to denote $a\le C b$ and  $a \sim b$ to denote $C^{-1}b\le a\le Cb$, where  $C$ is the universal constant  as defined above.

3) In the rest of the paper, we will use the notations
$$    \|\cdot\|_{L^p}=:\|\cdot\|_{L^p(I)},\q 1\le p\le \infty, \q {\rm and}\q \|\cdot\|_{H^k}=:\|\cdot\|_{H^k(I)},\q k>0.$$

\section{Lagrangian reformulation and main results}

To fix the boundary, we  define the {\it Lagrangian variable (trajectory)} $\eta^\e(x,t)$ by
\be\label{3}
\eta_t^\e(x,t)=u^\e(\eta^\e(x,t),t)\q{\rm for}\; t>0\q {\rm and}\q \eta^\e(x,0)=\eta_0^\e(x),\; x\in [0,\bar l],
\ee
where $\eta_0^\e: I\equiv(0, \bar l)\rightarrow(0,l_0)$ is a diffeomorphism (see \cite{OZ}) defined by
\be\label{eta0def}
\int_0^{\eta_0^\e(x)}\rho_0^\e(y)dy=\int_0^x\bar\rho (y)dy=M- g^{\frac{1}{\gamma-1}}\lt[\frac{\gamma-1}{\gamma}(\bar l-x)\rt]^\frac{\gamma}{\gamma-1},\q \eta_0^\e(\bar l)=l_0.
\ee
The above problem is equivalent to an initial value problem for $\eta_0^\e$:
\be \label{eta00}
\begin{split}
&\rho_0^\e(\eta_0^\e(x))\eta_{0x}^\e(x)=\bar\rho(x)=\lt[\frac{g(\gamma-1)}{\gamma}(\bar l-x)\rt]^{\frac{1}{\gamma-1}},\q \eta_0^\e(\bar l)=l_0.
\end{split}
\ee
Due to the assumptions \ef{initialdata}, \ef{oo} on $\rho_0^\e$ and $\bar\rho$, and that  the total mass of the steady solution is the same as that of $\rho_0^\e$, the diffeomorphism $\eta_0^\e$ is well defined.  Moreover, near the possible singular point $x=\bar l$, we have $ l_0-\eta_0^\e(x)=O(1) (\bar l-x)$.

Next, we set the {\it Lagrangian density} and {\it velocity}, respectively, by
\be\label{lag}
f^\e(x, t)=\rho^\e(\eta^\e(x, t), t) \ \ {\rm and} \ \  v^\e(x, t)=u^\e(\eta^\e(x,t), t).
\ee
Then the  system \ef{1} can be written as
\be\label{419}\lt\{ \begin{split}
& f_t^\e + f^\e \frac{v_x^\e}{\eta_x^\e}=0    & {\rm in}& \  \    I\times (0, \iy),\\
& f^\e v_t^\e+ \frac{1}{\epsilon^2}\frac{((f^\e)^\gamma)_x}{\eta_x^\e}=\frac{1}{\eta_x^\e}\lt(\frac{  v_x^\e}{\eta_x^\e}\rt)_x-\frac{1}{\epsilon^2}f^\e g \ \  &{\rm in}& \  \   I\times (0, \iy), \\
& f^\e>0 \ \  &{\rm in}& \  \   I\times (0, \iy), \\
& v^\e(0,t)=0, \  \    \lt(\frac{(f^\e)^\gamma}{\epsilon^2}-\frac{v_x^\e}{\eta_x^\e}\rt)(\bar l,t)=0,\ \ & {\rm for}&\ \  t>0,\\
& (f^\e, v^\e) =(\rho_0^\e(\eta_0^\e), u_0^\e(\eta_0^\e))  &  {\rm on}& \  \   I \times \{t=0\}.
\end{split}\rt.
\ee
Suppose that $\eta_x^\e>0$ for all $(x,t)\in I\times [0,\iy)$, which can be recovered later in the estimate \ef{11}. Then we   solve  $\eqref{419}_1$ to get
\be\label{5}
f^\e(x,t)=\frac{\bar\rho(x)}{\eta_x^\e(x,t)}, \q (x,t)\in I\times [0,\iy).
\ee
From \ef{oo}, without the loss of generality, we assume $\bar l=1$ and $g=1$, which gives $$\bar\rho(x)=\lt[\frac{\gamma-1}{\gamma}(1-x)\rt]^{\frac{1}{\gamma-1}}.$$

Formally, if $f^\epsilon\to \bar\rho$ or equivalently $\eta_x^\e\to 1$ as $\e\to 0$, we obtain from $\ef{419}_1$ that $\bar\rho v_x=0$, which is different from the one in the Eulerian coordinate. However, since  $\bar\rho$ only vanishes at $x=1$, we get $v_x=0$ $a.e.$. Due to $v(0,t)\equiv 0$, one derives that $v=0$ $a.e.$. Thus the limiting state is also $(\bar\rho,0)$.

Inserting \ef{5} into \eqref{419},  the system reads as
\be\label{6}\lt\{\begin{split}
& \bar\rho v_t^\e + \frac{1}{\epsilon^2}\lt( \frac{\bar\rho^\gamma}{(\eta_x^\e)^{\gamma}} \rt)_x = \lt(\frac{ {v_x^\e}}{\eta_x^\e}\rt)_x -\frac{1}{\epsilon^2}\bar\rho g \ \  & {\rm in} & \ \ I \times (0,\iy),\\
&v^\e(0,t)=0,\q v_x^\e(1,t)=0 & {\rm on} & \ \    (0,\iy),\\
&(\eta^\e, v^\e)(x,0)=(\eta_0^\e(x),  v_0^\e(x) )   & {\rm on} & \ \ I,
\end{split}\rt.
\ee
where $v_0^\e(x):=u_0^\e(\eta_0^\e(x))$ and $\eta^\e(x,t)$ satisfies \ef{3}.

\noindent
{\bf Notations.} For any $\epsilon\in (0,1]$ and $\alpha>0$, define the low-order energy, high-order energy and total energy, respectively, by
\be
\begin{split}
\mathfrak{E}_{L}^\e(t) := &\int_I [\bar\rho (v^\e)^2 +  \bar\rho^{1-\gamma+\alpha}(\eta_x^\e-1)^2 + \bar{\rho}^{\gamma-1+\alpha} (\eta_{xx}^\e)^2 ]dx(t)\\
&+ \int_I[\e^2((v_x^\e)^2+ (1+t)^{-\frac{\gamma-1}{\gamma}+\alpha}  (\eta_{xx}^\e)^2)+ \e^4\bar\rho (v_t^\e)^2 ]dx (t)\\
&+ \frac{1}{\e^2}\int_I     {\bar\rho^\gamma}  \lt(  ¡¡\eta_x^\e-1\rt)^2  dx (t)¡¡+ \|\eta_{x}^\e(\cdot,t)\|_{L^\infty(I)}+ \|(\eta_{x}^\e)^{-1}(\cdot,t)\|_{L^\infty(I)},
\end{split}
\ee
\be
\mathfrak{E}_{H}^\e(t):=\int_I (\e^8\bar\rho (v_{tt}^\e)^2 + \e^2 \bar{\rho}^{3\gamma-3+\alpha}(\eta_{xxx}^\e)^2 +\e^8 (1+t)^{-\frac{3\gamma-3+\alpha}{\gamma}} (\eta_{xxx}^\e)^2)dx(t),$$
\ee
\be\label{Et}
\mathfrak{E}^\e (t):=\mathfrak{E}_{H}^\e(t)+\mathfrak{E}_{L}^\e(t).$$
\ee

\begin{thm}\label{local} (Local Existence) Assume that $\gamma>1$ and fix $\e\in (0,1]$. Suppose that $(\eta_0^\e, v_0^\e)\in H^3(I)$, $(D_0^\epsilon)^{-1}\le\eta_{0x}^\e\le D_0^\epsilon$ for some constant $D_0^\e>0$, and there exist functions $h_1^\e\in H^2(I)$ and  $h_2^\e\in L^2(I),$  such that
\be\label{h1}
\bar\rho h_1^\e + \frac 1 {\e^2} \left(\frac{\bar\rho^\gamma}{(\eta_{0x}^\epsilon)^\gamma}\right)_x
=\left(\frac{v_{0x}^\e}{\eta_{0x}^\e}\right)_x-\frac{1}{\e^2}\bar\rho g,  \quad
\ee
 \be\label{h2}
\bar\rho h_2^\e - \frac \gamma {\e^2} \left(\frac{\bar\rho^\gamma v_{0x}^\e}{(\eta_{0x}^\epsilon)^{\gamma+1}}\right)_x
=\left(\frac{\eta_{0x}^\e h_{1x}^\e -(v_{0x}^\e)^2}{(\eta_{0x}^\e)^2}\right)_x,   \quad
\ee
and the following compatibility conditions hold:
\be\label{compatibility1}
v_0^\e(0)=0,\; v_{0x}^\e(1)=0,
\ee
\be\label{compatibility2}
 h_1^\e(0)=0,\; h_{1x}^\e(1)=0.
\ee
Then there exists a  constant $T^{*,\epsilon}>0$, such that \ef{6} admits a unique  solution $(\eta^\e,v^\e)$   satisfying
\be\label{regu}
  \begin{split}
  (\eta^\e,v^\e)\in L^\infty(0,T^{*,\epsilon};H^3(I)),\q v_t^\e\in L^\infty(0,T^{*,\epsilon}; W^{1,\infty} (I)),\\
  \sqrt{\bar\rho}v_{tt}^\e\in L^\infty(0,T^{*,\epsilon};L^2(I)),\q v_{tt}^\e\in L^2(0,T^{*,\epsilon};H^1(I)),\q
  \end{split}
  \ee
and $\underline{\lambda}^\e\le \eta_x^\e\le \bar\lambda^\e$ for some positive constants $\underline{\lambda}^\e$ and $\bar\lambda^\e$.
 \end{thm}
We will prove the above theorem in the final section by the method of finite difference scheme, which is also used in \cite{LiXY,LXZ,YZ}. Indeed, we have the regularity $v_t^\e\in L^\infty(0,T^{*,\epsilon}; H^2 (I))$ by \ef{regu} and $\ef{6}_1$.  We remark that although we only prove the existence of local classical solutions, however, the local strong solution satisfying \ef{lregu} can also be obtained in the same way.

Next, we state the main results of this paper.

\begin{thm}\label{thm2} (Uniform estimates and global existence for ``ill-prepared'' initial data) Suppose that $\e\in (0,1]$ is an arbitrary constant. Let $T$ be a given positive constant, and  $C,C_0,C_1$ are  generic positive constants independent of $\e\in (0,1]$ and $t\in [0,+\infty)$.

(i)
  Assume that $\gamma>1$,  $(\eta_0^\e, v_0^\e)\in H^2(I)$, \ef{compatibility1} and \ef{h1} are satisfied for $h_1^\e\in L^2(I)$, and $\mathfrak{E}_{L}^\e(0)\le C_0 $. Then there exists a unique global solution $(\eta^\e,v^\e)$ to \ef{6}  satisfying
    \be \label{lregu}
  \begin{split}
  &(\eta^\e,v^\e)\in L^\infty(0,T;H^2(I)),\; (\eta_x^\e)^{-1}\in L^\infty(0,T;L^\infty(I)),\\
  &\sqrt{\bar\rho}v_t^\e\in L^\infty(0,T;L^2(I)),\; v_t^\e\in L^2(0,T;H^1(I)),
  \end{split}
  \ee
  with the uniform-in-$\e$ estimates
  \be\label{e}
  \mathfrak{E}_{L}^\e(t)\le C\exp(C\bar Q^2)  (\mathfrak{E}_L^\e(0)+ (\mathfrak{E}_L^\e(0))^3), \quad for \quad t\ge 0, \; x\in I,
   \ee
   and
   \be\label{e1}
   \begin{split}
\frac{1}{\e^2}&\int_0^t\lt( \int_I  \bar\rho^{1+\alpha}    (\eta_x^\e -1)^2dx +  (1+s)^{\frac{\gamma-1}{\gamma}-\alpha } \int_I \bar\rho^\ga  (\eta_x^\e -1)^2dx\rt)ds\\
&
+\int_0^t (1+s)^{\frac{2\gamma-1}{\gamma}-\alpha}\int_I  ((v^\e)^2 + (v_x^\e)^2   )dx   ds\\
 &+\frac{1}{\e^2}\int_0^t\lt(\int_I\bar{\rho}^{2\gamma-1+\alpha}(\eta_{xx}^\e)^2dx+
 (1+s)\int_I\bar{\rho}^{3\gamma-1+\alpha}(\eta_{xx}^\e)^2dx\rt)ds \\
 &\le C\exp(C\bar Q^2)  (\mathfrak{E}_L^\e(0)+ (\mathfrak{E}_L^\e(0))^3).
\end{split}
   \ee
Moreover, the following uniform decay estimates are satisfied:
\be\label{decay1}
\begin{split}
(1&+t)^{\frac{1}{\gamma}[\min(\gamma-1,\frac{\gamma}{2})+\gamma-1]-\alpha}
\|\eta^\e(\cdot,t)-x\|_{L^\infty}^2 \\
&+(1+t  )^{\frac{1}{\gamma}   [\min(2\gamma-2,\gamma)+\gamma-1]-\alpha}  \|\eta^\e(\cdot,t)-x\|_{L^2}^2 \\
&+ (1  +t)^{\frac{\gamma-1}{\gamma}-\alpha }\|\eta_x^\e(\cdot,t)  -1\|_{L^2 }^2
  +     (1+t)^{\frac{2\gamma-1}{\gamma}-\alpha} \|\lt(\e^2\sqrt{\bar\rho} v_t^\e, \e^2 v_x^\e \rt)(\cdot,t) \|_{L^2 }^2\\
 &+(1+t)^{\frac{2\gamma-1}{\gamma}-\varsigma}\e^4\|v_x(\cdot,t)\|_{L^{\infty}}^2+
 (1+t)^{ \frac{1}{\gamma} \min{(1,\gamma-1) }-\alpha}\|\e^2 v_{xx}^\e(\cdot,t)\|_{L^2}^2  \\
&\quad\le C \exp (C\bar Q^2) (\mathfrak{E}_L^\e(0)+ (\mathfrak{E}_L^\e(0))^3  ),
  \end{split}
\ee
for any constant $\alpha>0$.

  (ii) Suppose in addition that $\gamma\in (1,3/2)$, $(\eta_0^\e, v_0^\e)\in H^3(I)$, \ef{h1}-\ef{compatibility2} hold for $h_1^\e\in H^2(I)$ and  $h_2^\e\in L^2(I),$ and   $\mathfrak{E}_{H}^\e(0)\le C_1.$ Then the solution $(\eta^\e,v^\e)$ obtained in (i) is a classical solution to \ef{6}, satisfies that
  \be\label{hregu}
  \begin{split}
  &(\eta^\e,v^\e)\in L^\infty(0,T;H^3(I)),\; (\eta_x^\e)^{-1}\in L^\infty(0,T;L^\infty(I)),\; v_t^\e\in L^\infty(0,T;W^{1,\infty}),\\
  &\sqrt{\bar\rho}v_{tt}^\e\in L^\infty(0,T;L^2(I)),\q v_{tt}^\e\in L^2(0,T;H^1(I)),
  \end{split}
  \ee
  with the uniform-in-$\e$ bounds
  \be\label{e2}
  \mathfrak{E}^\e(t)+\int_0^t\int_I\bar{\rho}^{4\gamma-3+\alpha}(\eta_{xxx}^\e)^2dxdt\le C\exp(C\bar Q^2)  (\mathfrak{E}^\e(0)+ (\mathfrak{E}^\e(0))^{10}),\quad t\ge 0,
  \ee
  and the following decay estimates:
  \be\label{decay2}
\begin{split}
(1&+t)^{\frac{2\gamma-1}{\gamma}-\alpha}\|\e^4\sqrt{\bar{\rho}}v_{tt}^\e(\cdot,t)\|_{L^2}^2
 +(1+t)^{\frac{2\gamma-1}{\gamma}-\alpha}\|\e^4v_{tx}^\e(\cdot,t)\|_{L^\infty}^2 \\
 &\le C \exp(C\bar Q^2)  (\mathfrak{E}^\e(0)+ (\mathfrak{E}^\e(0))^{10}),
\end{split}
\ee
for any constant $\alpha>0$.

\end{thm}

We show this theorem by the local existence result in Theorem \ref{local} and the uniform estimates established in Section 4.

\begin{rmk} Although the uniform estimates in \ef{e2} and \ef{decay2} are not used directly in the convergence $\epsilon\to 0$ in Theorem 2.6 below, however, it gives the global existence and decay rates of classical solutions of (2.7) (for fixed $\e$), moreover, provides us some information about the variation of the pressure. Indeed, from the uniform weighted estimates, we can obtain the uniform estimates and the convergent rates for the variation of the pressure and its spatial derivatives as follows:
\bee
\begin{split}
&\|p(f^\e)-p(\bar\rho)\|_{L^2}^2(t)\lesssim\int_I\bar\rho^\gamma(\eta_x^\e-1)^2dx(t)\lesssim \e^2,\\
&\|(p(f^\e)-p(\bar\rho))_x\|_{L^2(I\times(0,t))}\lesssim  \int_0^t\int_I (\bar\rho^2(\eta_x^\e-1)^2 +\bar\rho^{2\gamma-2}(\eta_{xx}^\e)^2)dxdt\lesssim\e^2\\
&\|(p(f^\e)-p(\bar\rho))_{xx}\|_{L^2(I\times(0,t))}\\
 &\quad\lesssim \int_0^t  \int_I  \lt[\bar\rho^{1+\alpha}    (\eta_x^\e -1)^2 +  \bar\rho^{2\gamma-1+\alpha}    (\eta_{xx}^\e )^2 +  \bar\rho^{4\gamma-3+\alpha}    (\eta_{xxx}^\e )^2\rt]dx  ds\lesssim C,
\end{split}
\eee
where the condition $\gamma<\frac 3 2$ is required in the last inequality.
\end{rmk}

\begin{rmk}
The assumption $\gamma\in (1,3/2)$ is physical since it includes the case of the air where $\gamma\approx 1.4$.
\end{rmk}

 \begin{rmk}\label{classical} There are no uniform-in-$\epsilon$ estimates for the high norms of the regularity \ef{hregu}. However, for each fixed $\epsilon$, the regularity in \ef{hregu} and  Lemma \ref{compact} give that
\bee
\begin{split}
& v^\e\in C([0,T],C^{2+\alpha}(I) ) ,\;  \bar\rho v_t^\e\in C([0,T], C^{\alpha} (I) ),\;\eta^\e\in C^1([0,T],C^{2+\alpha}(I) ),
\end{split}
\eee
for any $\alpha\in (0,1/2)$, thus the global solution $(\eta^\e,v^\e)$ obtained in Theorem \ref{thm2} is a classical solution to \ef{6}. Therefore $(\rho^\e,u^\e)$ is a classical solution to \ef{1} since $\eta_x^\e$ is bounded from below and above by \ef{e}, i.e., the mapping $(\eta,t)\mapsto (x,t)$ is a bijection.
\end{rmk}

\begin{thm}\label{thm3} (Singular limit for ``ill-prepared'' data) Suppose that the assumptions in Theorem \ref{thm2} (i) are satisfied and  $(\eta^\epsilon, v^\epsilon)(x,t)$ is a strong solution obtained in Theorem \ref{thm2} (i) for each $\e$. Then  we obtain
\be\label{conv1}
  \eta^\epsilon\to x,\;\eta_x^\epsilon \to 1\; {in} \; C([0,T],L^2(I)),\q   v^\epsilon \rightharpoonup 0,\;v^\epsilon_x \rightharpoonup 0  \;   {in}\;  L^2(0,T;L^2(I)),\q  \epsilon\to 0 .
  \ee
Thus the solution $(\rho^\epsilon,u^\epsilon)$ of \ef{1} converges in the following sense: \be\label{conv2}
(\rho^\epsilon(\eta^\epsilon(x,t),t),u^\epsilon(\eta^\epsilon(x,t),t))\to(\bar\rho(x),0)\q a.e.  \;(x,t)\in I\times[0,T],\q\epsilon\to 0.
 \ee
 Moreover, the free boundary
 \be\label{conv3}
 \eta^\epsilon(1,t)=\Gamma^\e(t)\to 1 \q  {\rm for}\q  a.e. \;  t\in[0,T],\q \epsilon\to 0.
 \ee
\end{thm}

The proof will be given in Section 5.

\begin{rmk}
In Theorem \ref{thm3}, the initial data are ``ill-prepared'' since  ``$\bar\rho v_t^\e(x,0)$''$\neq O(1)$ and $\rho_0^\e-\bar\rho=O(\e)$. The convergence \ef{conv2} is a direct consequence of \ef{conv1} since $\rho^\epsilon(\eta^\epsilon(x,t),t)=f^\e(x,t)=\bar\rho(x)(\eta_x^\e(x,t))^{-1}$ and $u^\epsilon(\eta^\epsilon(x,t),t)=v^\e(x,t)$. The convergence of $\{(\rho^\epsilon,u^\epsilon)\}_{\e>0}$ is not in the classical sense, since the domain of the limiting state $(\bar\rho, 0)$ is different from those of $\{(\rho^\epsilon,u^\epsilon)\}_{\e>0}$  which are time-dependent.
\end{rmk}

\noindent{\bf Notation.}
For any $\epsilon\in (0,1]$ and $\alpha>0$, we define the ``well-prepared'' low-order energy by
\bee
\begin{split}
\tilde{\mathfrak{E}}_{L}^\e(t) := &\frac{1}{\e^4}\int_I [\bar\rho (v^\e)^2 +   \bar\rho^{1-\gamma+\alpha}(\eta_x^\e-1)^2]dx (t)\\
&+\int_I (\frac{1}{\e^2}(v_x^\e)^2+  \bar\rho  (v_t^\e)^2)dx(t)+ \frac{1}{\e^6}\int_I     {\bar\rho ^\gamma}  \lt(  \frac{1}{\gamma-1}{(\eta_x^\e)^{1-\gamma}}+ \eta_x^\e-\frac{\gamma}{\gamma-1}\rt)  dx (t).
\end{split}
\eee

\begin{coro}
(Uniform estimates and strong convergence for ``well-prepared'' initial data)

 \noindent Assume that $\gamma>1$, $(\eta_0^\e, v_0^\e)\in H^2(I)$, \ef{compatibility1} and \ef{h1} are satisfied for $h_1^\e\in L^2(I)$,  $\tilde{\mathfrak{E}}_{L}^\e(0)\le C_0 $ and $C_0^{-1}\le \eta_{0x}^\e\le C_0$. Then the solution $(\eta^\e,v^\e)$   satisfies
  \bee\label{e}
  \tilde{\mathfrak{E}}_{L}^\e(t)+\frac 1 {\e^4}\int_0^t (1+s)^{\frac{2\gamma-1}{\gamma}-\alpha}\int_I  [(v^\e)^2 + (v_x^\e)^2   ]dx   ds\le C ,\quad C^{-1}\le\eta_x^\e(x,t)\le C \quad for \; t\ge 0, \; x\in I.
   \eee
Moreover, the following decay estimates are satisfied:
\bee
\begin{split}
  (1&+t)^{\frac{2\gamma-1}{\gamma}-\alpha}\|  v_x^\e(\cdot,t)\|_{L^{\infty}}^2
 +     (1+t)^{\frac{2\gamma-1}{\gamma}-\alpha} \|\lt( \sqrt{\bar\rho} v_t^\e,   v_x^\e \rt)(\cdot,t) \|_{L^2 }^2\le C,
  \end{split}
\eee
for any constant $\alpha>0$.
Therefore for any $T>0$, we obtain that $v^\epsilon\to 0$ in $L^2(0,T;H^1(I))$, $\eta^\epsilon\to x,\;\eta_x^\epsilon\to 1$ in $C([0,T],L^2(I))$ as $\e\to 0$,  with the rate $\e^2$.
\end{coro}
The proof of this corollary is exactly the same as Theorem \ref{thm2} (i), thus we omit the details.

\section{Preliminaries}
In this sub-section, we present some weighted inequalities.  The following general version of the Hardy inequality whose proof can be found in \cite{KM}
will  be used often in this paper.
\begin{lem} (Hardy's inequality)\label{hardy} Let $I=(0,1)$. For any given real number $\beta>-(\ga-1)$,
\be\label{hard}
\int_I \bar\rho^{\beta} w^2 dx
 \le C_1 \int_I \bar\rho^{\beta+2(\ga-1)}(w^2 + w_x^2)dx,
 \ee
provided that  $\int_I \bar\rho^{\beta+2(\ga-1)}(w^2 + w_x^2)dx<\iy$,  where $C_1$ is a generic constant independent of $w$.
\end{lem}

\begin{lem}\label{hardy2} (Refined Hardy's inequality) {\rm(See \cite{OZ})}  Let $I=(0,1)$ and $\beta>-(\ga-1)$ be a given real number and $w(x)$ be a function  satisfying $w(0)=0$. Then it holds
\be\label{lem33}
\int_I \bar\rho^\beta w^2 dx\le C_2 \int_I  \bar\rho^{\beta + 2(\ga-1)} w_x^2 dx,
\ee
where $C_2$ is a generic constant independent of $w$.
\end{lem}
To show the regularity of solutions, we need the following compactness theorem.
\begin{lem}\label{compact}{\rm (See \cite{simon})} Suppose that $X\subset Z\subset Y$ are Banach spaces and the embedding $X\subset Z$ is compact. Let $E:=\{u(x,t)|u\in L^\infty(0,T;X),u_t\in L^r(0,T;Y)\}$ for $r>1$. Then the embedding $E\subset C([0,T], Z)$ is compact.
\end{lem}

Next, we prove a  weighted inequality by using previous lemma, with additional assumptions on $\bar\rho(x)$. This inequality can be used to  prove the uniqueness.
\begin{lem}\label{lem21} Assume that $\bar\rho(x)$ is a continuous function in $I=[0,1]$, satisfying $\bar\rho(x)> 0 $ for $x\in [0,1)$ and $\bar\rho(1)=0$. Then for any  $\epsilon>0$, there exists a positive constant $C=C(\epsilon)$  depending on $\epsilon$, but not on $u$, such that
\be\label{B9}
\|u\|_{L^2(I)}^2\le \epsilon\|u_x\|_{L^2(I)}^2 +    C  (\epsilon) \int_I\bar\rho u^2 dx , \q \;u\in H^1(I).
\ee
\end{lem}
{\bf Proof.}  We prove this inequality by dividing
\be\label{B10}
 \int_I  u^2dx =   \int_{1-\delta}^{1}  u^2 dx  +  \int_0^{ 1-\delta} u^2 dx \equiv K_1+K_2,
\ee
for any $\delta\in (0,1)$. Here
\be\label{B11}
\begin{split}
K_1 &\lesssim  \int_{1-\delta}^{1} \| u \|_{L^\infty(I)}^2 dx \\
&\le  \delta  \lt[\| u_x \|_{L^1(I)}^2 + \lt(\int_I\bar\rho |u| dx\rt)^2  \rt] \\
&\le  C_1\delta \lt(  \int_I   u_x^2 dx   +   \int_I \bar\rho  u^2 dx  \rt).
\end{split}\ee
 On the other hand, since the continuous function $\bar\rho(x)>0$ for $x\in [0, 1)$, there exists a positive constant $ {\rho}_\delta$, such that
$$\bar\rho(x)\ge  {\rho}_\delta>0, \q x\in [0,1-\delta].$$
Thus we get
\be\label{B12}
\begin{split}
K_2  \le    \int_0^{1-\delta} \frac{\bar\rho}{ {\rho}_\delta} u^2 dx
 \le   {\rho}_\delta^{-1}    \int_I   \bar\rho  u^2 dx .
\end{split}\ee
Submitting \ef{B11} and \ef{B12} into \ef{B10}, we obtain \ef{B9}, where $\epsilon =C_1\delta$ and $C(\epsilon)=C_1\delta + {\rho}_\delta^{-1} .$

\hfill$\Box$

\section{Uniform estimates for ``ill-prepared'' initial data}

Suppose that $\gamma>1$ and let $\epsilon$ be an arbitrary constant in $(0,1]$. We also assume that $(\eta^\epsilon,v^\epsilon)(x,t)$ is a sufficiently smooth solution to \ef{6} in $I\times[0,T]$, where $T>0$ is an arbitrary constant.  Moreover, $\eta^\epsilon(x,t)$ satisfies  \ef{3} with
\be\label{assum}
\eta_x^\epsilon(x,t)>0,\q {\rm for}\q (x,t)\in I\times[0,T],
\ee
which will be recovered in Lemma \ref{lem2}. For simplicity, we  drop the superscript $\epsilon$ from now on.

\subsection{Uniform lower-order estimates and global strong solutions}
\begin{lem}\label{lem1} {\rm (Basic energy)}
Assume \ef{assum} holds and let $$E  (t):=\int_I\lt[\frac 12  \bar\rho v^2 +   \frac{\bar\rho^\gamma}{\e^2} \lt(  \frac{1}{\gamma-1}{\eta_x^{1-\gamma}}+ \eta_x-\frac{\gamma}{\gamma-1}\rt)\rt] dx (t).$$ Then for any $(x,t)\in I\times[0,T],$ we have
\be\label{7A}
\int_I\bar\rho v^2 dx(t)+  \int_0^t\int_I \eta_x^{-1} v_x^2 dxdt\le E_0 ,
\ee
where $E_0 :=\int_I\lt[\frac 12  \bar\rho v_0^2 +   \frac{\bar\rho^\gamma}{\e^2} \lt(  \frac{1}{\gamma-1}{\eta_{0x}^{1-\gamma}}+ \eta_{0x}-\frac{\gamma}{\gamma-1}\rt)\rt] dx\le \mathfrak{e}_0<+\infty,$ with $\mathfrak{e}_0$ being a positive constant independent of $\e\in(0,1]$ and $t\in[0,+\infty)$.
\end{lem}
{\bf Proof.} Using \eqref{station},   we rewrite $\ef{6}_1$ as
\be\label{8}
\bar\rho v_t + \frac{1}{\e^2}\lt[ \bar\rho^\gamma( \eta_x^{-\gamma}-1)\rt]_x = \lt(  \eta_x^{-1} {v_x}\rt)_x .
\ee
Observe that
$$( \eta_x^{-\gamma}-1) v_x=-\lt(   \frac{1}{\gamma-1}{\eta_x^{1-\gamma}}+\eta_x -\frac{\gamma}{\gamma-1}  \rt)_t.$$
We multiply \ef{8} by $v$  and integrate the resulting equality to get
\be\label{9}\begin{split}
  \frac{d}{dt} E(t) +    \int_I   \eta_x^{-1} {v_x^2}dx  =0,
\end{split}\ee
where  the boundary conditions $v(0,t)=0$ and $\bar\rho(1)=0$ are used. Note that by the Taylor expansion, we have
 \be\label{Taylor}
 \begin{split}
 \frac{1}{\gamma-1}{\eta_x^{1-\gamma}}+ \eta_x-\frac{\gamma}{\gamma-1}&=\frac{1}{\gamma-1}[{\eta_x^{1-\gamma}} -1 -(1-\gamma)(\eta_x-1)]\\
 &= \gamma\xi^{-\gamma-1}(\eta_x-1)^2\ge C(\e,T)(\eta_x-1)^2\ge 0,
 \end{split}
 \ee
 where $\xi$ is a function between $\eta_x$ and $1$ satisfying $0<\xi<1+\max_{x,t}\eta_x(x,t)\le C(\e,T)$.
 Thus \ef{7A} follows. \hfill$\Box$
\vskip 3mm
In the following, we show a key lemma of this paper. It indeed implies that the size of the free interval is bounded from below and above, where the bounds are independent of $\e\in(0,1]$ and $t\in(0,+\infty)$.
\begin{lem}\label{lem2}
Suppose that  \ef{assum} is satisfied  and $D_0^{-1}\le \eta_{0x}\le D_0$. Then
\be\label{11}
0< \bar Q^{-1}  \le \eta_x (x ,t) \le   \bar Q  , \q x\in I,\;t\in [0,T],
\ee
where  $\bar Q:=D_0^2\exp\lt\{4\sqrt{M\mathfrak{e}_0} \rt\}$ and  the constant $\mathfrak{e}_0$ is defined in Lemma 4.1.
 \end{lem}
{\bf Proof.} Integrating \ef{8} with respect to the spatial variable in $(x,1)$, we find
\be \label{13}
\begin{split}
 \lt(\log\eta_x \rt)_t& =- \int_{x}^1 \bar\rho v_t dy +\frac{\bar\rho^\gamma}{\e^2}\lt(\eta_x^{-\gamma}-1\rt)\\
 &=-\frac{d}{dt} \int_{x}^1 \bar\rho v  dy +\frac{\bar\rho^\gamma}{\e^2}\lt( \eta_x^{-\gamma}-1\rt).
\end{split}\ee
Next, integrate \ef{13}   from $0$ to $t$ to get
\be\label{14}
\begin{split}
\eta_x&=\eta_{0x}\exp\lt\{-\int_x^1 \bar\rho(y) v(y,t)dy\Big|_0^t\rt\}\exp\left\{\frac{\bar\rho^\gamma}{\e^2} \int_0^t (\eta_x^{-\gamma}-1) dt\rt\}\\
&\equiv \eta_{0x}\mathcal{Q}(x,t)\mathcal{R}(x,t).
\end{split}
\ee
 Then we calculate the highly oscillatory term $\mathcal{R}(x,t)$ by
 \be\label{15}
\begin{split}
\frac{\partial \mathcal{R}}{\partial t}&=\frac{\bar\rho^\gamma}{\e^2} (\eta_x^{-\gamma}-1)\mathcal{R} = -\frac{\bar\rho^\gamma}{\e^2} \mathcal{R} +\frac{\bar\rho^\gamma}{\e^2} \eta_{0x}^{-\gamma}\mathcal{Q}^{-\gamma} \mathcal{R}^{1-\gamma}.
\end{split}
\ee
Next, we introduce $\mathcal{Z}(x,t):= \mathcal{R}^\gamma(x,t)$. Then \ef{15} reads
\bee
\frac{\partial \mathcal{Z}}{\partial t}+\frac{\gamma}{\e^2}\bar\rho^\gamma\mathcal{Z} =\frac{\gamma}{\e^2} \bar\rho^\gamma \eta_{0x}^{-\gamma}\mathcal{Q}^{-\gamma}.
\eee
Solving this ordinary differential equation (for fixed spatial variable), we find
 \be\label{16}
\begin{split}
    \mathcal{Z}^\frac{1}{\gamma}=\mathcal{R}&= \exp\lt\{-\frac{\bar\rho^\gamma}{\e^2} t\rt\}\lt[1 +\frac{\gamma}{\e^2} \int_0^t\bar\rho^\gamma \eta_{0x}^{-\gamma}\mathcal{Q}^{-\gamma} \exp\lt\{\frac{\gamma}{\e^2}\bar\rho^\gamma s\rt\}ds\rt]^\frac{1}{\gamma}.
\end{split}
\ee
 Note that
\be\label{17A}
\begin{split}
\mathcal{Q}(x,t)&\le\exp\lt\{\lt(\int_I\bar\rho dx\rt)^{\frac 12}\lt[\lt(\int_I\bar\rho v^2dx\rt)^{\frac 12}+ \lt(\int_I\bar\rho v_0^2 dx\rt)^{\frac 12}\rt]\rt\}\\
&\le \exp\lt\{2\sqrt{M\mathfrak{e}_0} \rt\},
\end{split}
\ee
and similarly
\be\label{17B}
\begin{split}
\mathcal{Q}(x,t)&\ge \exp\lt\{-2\sqrt{M\mathfrak{e}_0} \rt\}.
\end{split}
\ee
Plug \ef{16}, \ef{17A} and \ef{17B} into \ef{14} to get
\bee
\begin{split}
\eta_x& \le  \eta_{0x} \exp\lt\{2\sqrt{M\mathfrak{e}_0} \rt\} \exp\lt\{-\frac{\bar\rho^\gamma}{\e^2} t\rt\}\\
&\qquad \times\lt[1 + \eta_{0x}^{-\gamma}\exp\lt\{2\gamma\sqrt{M\mathfrak{e}_0} \rt\} \int_0^t \frac{\gamma}{\e^2}\bar\rho^\gamma  \exp\lt\{\frac{\gamma}{\e^2}\bar\rho^\gamma s\rt\}ds\rt]^\frac{1}{\gamma}\\
   &\le   \bar Q^{\frac 12}\lt[\exp\lt\{-\frac{\gamma}{\e^2}\bar\rho^\gamma t\rt\}+ \bar Q^{\frac \gamma 2}(1-\exp\lt\{-\frac{\gamma}{\e^2}\bar\rho^\gamma t\rt\})\rt]^\frac{1}{\gamma} \le  \bar Q  ,
\end{split}
\eee
and
\bee
\begin{split}
   \eta_x &\ge  \eta_{0x} \exp\lt\{-2\sqrt{M\mathfrak{e}_0} \rt\} \exp\lt\{-\frac{\bar\rho^\gamma}{\e^2} t\rt\}\\
&\qquad \times\lt[1 + \eta_{0x}^{-\gamma}\exp\lt\{-2\gamma\sqrt{M\mathfrak{e}_0} \rt\} \int_0^t \frac{\gamma}{\e^2}\bar\rho^\gamma  \exp\lt\{\frac{\gamma}{\e^2}\bar\rho^\gamma s\rt\}ds\rt]^\frac{1}{\gamma}\\
   &\ge \bar Q^{-\frac 12}\lt[\exp\lt\{-\frac{\gamma}{\e^2}\bar\rho^\gamma t\rt\}+ \bar Q^{-\frac \gamma 2}(1-\exp\lt\{-\frac{\gamma}{\e^2}\bar\rho^\gamma t\rt\})\rt]^\frac{1}{\gamma} \ge    \bar Q^{-1} .
\end{split}
\eee
This  finishes the proof.
 \hfill$\Box$

\vskip 0.5cm
\begin{rmk}
The estimate \ef{11} recovers the assumption \ef{assum}. Thus in what follows, we don't need to suppose \ef{assum} in the energy estimates.
\end{rmk}

\begin{lem}\label{lem3} For $t\ge 0$, we have
\be\label{21}
\begin{split}
    \int_I  &   (\eta_x  -1)^2  dx(t)   + (1+t) \int_I\lt[ \bar\rho v^2 + \frac{\bar\rho^\gamma}{\e^2} (\eta_x-1)^2\rt] dx(t)\\
    &+ \int_0^t  \lt[\int_I \frac{\bar\rho^{ \gamma}}{\e^2}(\eta_x -1)^2dx  + (1+s)\int_I (v^2 + v_x^2 )dx \rt]  ds   \lesssim   \bar Q^{4 \gamma+6 }\mathfrak{e}_0  .
 \end{split}
\ee
\end{lem}
{\bf Proof.} Similar to \ef{Taylor}, using Taylor's expansion  and \ef{11}, one finds
\be\label{22A}
\begin{split}
  \bar Q^{- (\gamma+1)} (\eta_x-1)^2&\lesssim\eta_x +\frac{1}{\gamma-1}\eta_x^{1-\gamma}-\frac{\gamma}{\gamma-1}\\
  &=\gamma \xi^{-\gamma-1}(\eta_x-1)^2\lesssim \bar Q^{ \gamma+1 } (\eta_x-1)^2,
  \end{split}
\ee
where $\xi$ is between $\eta_x$ and $1$. Then \ef{9}  gives
\be\label{22C}
E(t)+ \int_0^t  \int_I   \eta_x^{-1} {v_x^2}dxdt \lesssim  \bar Q^{\gamma+1} \int_I \lt(\bar\rho v_0^2+ \frac 1 {\e^2} \bar\rho^\gamma(\eta_{0x}  -1)^2\rt) dx .
\ee
Note that $\eta_t=v$,  and thus rewrite \ef{8} as
 \be\label{19}
 -(\log\eta_x )_{xt} + \lt[\frac{\bar\rho^\gamma}{\e^2}(\eta_x^{-\gamma}-1)\rt]_x+\bar\rho v_t=0.
 \ee
Multiplying \ef{19} by a new multiplier
$$ \mathcal{G}(x,t;\beta):= -\int_0^x \bar\rho^{-\beta}(y) \log\eta_y (y,t) dy $$
 with the constant $\beta<  \gamma-1$ and integrating the resulting equality over $I$ give that
 \be\label{20}
\begin{split}
 \frac{1}{2 } \frac{d}{dt} &\int_I \bar\rho^{ -\beta}(\log\eta_x )^2 dx +
\frac{1}{\e^2}\int_I \bar\rho^{\gamma-\beta}(1-\eta_x^{-\gamma})\log\eta_x dx\\
 &= \int_I\bar\rho v_t\lt(\int_0^x \bar\rho^{-\beta} \log\eta_y  dy \rt)dx.
 \end{split}
\ee
We remark that this multiplier is motivated from   \cite{LXZ}.

Choose $\beta=0$ in \ef{20} and integrate with respect to the time variable to get
\bee
\begin{split}
   \int_I  (\log & \eta_x )^2 dx (t) + \frac{1}{\e^2}\int_0^t  \int_I \bar\rho^{\gamma }(1-\eta_x^{-\gamma})\log\eta_x dx dt\\
 \lesssim & \int_I  (\log\eta_{0x} )^2 dx+ \lt| \int_0^t \int_I \bar\rho v_t \lt(\int_0^x \log\eta_y  dy \rt)dxdt\rt|\\
 = & \int_I  (\log\eta_{0x} )^2 dx+  \lt|\int_I \bar\rho v\lt(\int_0^x \log\eta_y  dy \rt)dx\Big|_0^t-\int_0^t\int_I \bar\rho v\lt( \int_0^x \eta_y^{-1} v_y dy\rt) dx dt\rt|\\
 \le & \delta \int_I (\log\eta_x)^2 dx(t) + C \lt( \int_I (\log\eta_{0x})^2 dx + \int_I \bar\rho v_0^2 dx\rt)\\
 & + C(\delta)  \int_I \bar\rho v^2 dx(t) + C\|\eta_x^{-1}\|_{L^\infty_{x,t}}^\frac 1 2\|\eta_x \|_{L^\infty_{x,t}}^\frac 1 2\int_0^t\int_I \eta_x^{-1}v_x^2 dxdt,
\end{split}
\eee
since
\bee
\begin{split}
  \Big|\int_0^t\int_I & \bar\rho v\lt( \int_0^x \eta_y^{-1} v_y dy\rt) dx dt\Big|\\
 \lesssim &  \int_0^t\lt(\int_I v^2 dx\rt)^{\frac 1 2}\lt(\int_I \eta_x^{-1}v_x^2 dx\rt)^{\frac 1 2}\|\eta_x^{-1}\|_{L^\infty}^\frac{1}{2}dt\\
 \lesssim &  \int_0^t\lt(\int_I v_x^2 dx\rt)^{\frac 1 2}\lt(\int_I \eta_x^{-1}v_x^2 dx\rt)^{\frac 1 2}\|\eta_x^{-1}\|_{L^\infty}^\frac{1}{2}dt\\
 \lesssim &  \int_0^t \int_I \eta_x^{-1}v_x^2 dx dt\|\eta_x^{-1}\|_{L^\infty_{x,t}}^\frac 1 2\|\eta_x \|_{L^\infty_{x,t}}^\frac 1 2.
\end{split}
\eee
By Taylor's expansions
$$\log y=\log 1 + \xi_1^{-1} (y-1),\quad y^{-\gamma}=1-\gamma \xi_2^{-\gamma-1}(y-1),$$
where $\xi_1,\xi_2$ are between $y$ and 1,
 and the uniform bounds of $\eta_x$ in \ef{11}, we find
\be
\begin{split}
&  \bar Q^{-2}(\eta_x-1)^2\lesssim (\log \eta_x)^2\lesssim \bar Q^2(\eta_x-1)^2,\\
&   \bar Q^{- (\gamma+2)}(\eta_x-1)^2\lesssim(1-\eta_x^{-\gamma})\log\eta_x\lesssim \bar Q^{ \gamma+2 }(\eta_x-1)^2.
\end{split}
\ee
Then, by use of \ef{22C} and Lemma \ref{lem2}, we obtain
\be\label{22}
\begin{split}
    \int_I  &   (\eta_x  -1)^2  dx (t)
    + \frac{1}{\e^2}\int_0^t  \int_I \bar\rho^{ \gamma}(\eta_x -1)^2dx ds   \lesssim\bar Q^{2\gamma+4}   \mathfrak{e}_0  .
 \end{split}
\ee
Next, we integrate $(1+t)\cdot\ef{9}$ with respect to the time variable to derive
\be\label{23A}
\begin{split}
  (1&+t) E(t) + \int_0^t(1+s)\int_I  \eta_x^{-1}v_x^2 dxds= E_0   + \int_0^tE(s)ds \\
& \lesssim  \bar Q^{  \gamma+1  } \lt(   \mathfrak{e}_0   + \int_0^t\int_I\lt[\bar\rho^{2\gamma-1}   v_x^2  + \frac{1}{\e^2}\bar\rho^\gamma(\eta_x-1)^2\rt]dxds \rt)\\
&
 \lesssim  \bar Q^{3 \gamma+5 }  \mathfrak{e}_0 ,
 \end{split}
\ee
by employing    \ef{11}, \ef{22A}, \ef{22C}, \ef{22} and Lemma \ref{lem33}.
Therefore, we obtain \ef{21} by use of \ef{11}, \ef{22}, \ef{23A} and Lemma \ref{lem33}.  \hfill$\Box$

\begin{lem} {\rm (faster decay)}\label{lem5} It holds that, for
any $t \ge 0$,
\be\label{30B}
\begin{split}
\int_I     \bar\rho^{1-\gamma +\varsigma} (\eta_x  -1)^2 dx (t)
+ \frac{1}{\e^2}\int_0^t \int_I  \bar\rho^{1+\varsigma}    (\eta_x -1)^2dx ds\lesssim  \exp(C\bar Q^2)  \mathfrak{E}_L(0),
 \end{split}
\ee
\be\label{30}
\begin{split}
         \int_I \Big[  (1& +t)^{\frac{\gamma-1}{\gamma}-\varsigma }  + (1+t)^{\frac{2\gamma-1}{\gamma}-\varsigma } \frac{\bar\rho^\ga}{\e^2}\Big](\eta_x  -1)^2 dx (t) +\e^2 (1+t)^{\frac{2\gamma-1}{\gamma}-\varsigma}\int_I v_x^2 dx (t)\\
& + \e^2\int_0^t (1+s)^{\frac{2\gamma-1}{\gamma}-\varsigma}\int_I\bar\rho v_s^2 dxds+\frac{1}{\e^2}\int_0^t   (1+s)^{\frac{\gamma-1}{\gamma}-\varsigma } \int_I \bar\rho^\ga  (\eta_x -1)^2dxds
\\
& +\int_0^t (1+s)^{\frac{2\gamma-1}{\gamma}-\varsigma}\int_I  (v^2 + v_x^2   )dx   ds \lesssim \exp(C\bar Q^2)  \mathfrak{E}_L(0)   ,
 \end{split}
\ee
\be\label{36}
\begin{split}
    (1+t& )^{\frac{1}{\gamma}   [\min(2\gamma-2,\gamma)+\gamma-1]-\varsigma}  \|\eta-x\|_{L^2}^2 +(1+t)^{\frac{1}{\gamma}[\min(\gamma-1,\frac{\gamma}{2})+\gamma-1]-\varsigma}\|\eta-x\|_{L^\infty}^2\\
  &\lesssim \exp (C\bar Q^2) \mathfrak{E}_L(0),
\end{split}
\ee
where $\varsigma\in (0, \frac{\ga-1}{\gamma})$ is an arbitrary constant.
 \end{lem}
\noindent{\bf Proof}. Similar as \ef{22}, for $\beta=\gamma-1-\varsigma$, $\varsigma>0$,  in \ef{20}, we have
\bee\label{}
\begin{split}
   \int_I  \bar\rho^{  1-\gamma+\varsigma}&(\eta_x  -1)^2 dx (t) +\frac{1}{\e^2}
\int_0^t\int_I \bar\rho^{1+\varsigma}(\eta_x -1)^2dxdt\\
 \lesssim & \bar Q^{\gamma+4}   \int_I \bar\rho^{1-\gamma+\varsigma} (\eta_{0x}  -1)^2  dx + \bar Q^{\gamma+2} \lt|\int_0^t \int_I\bar\rho v_t \int_x^1\bar\rho^{ 1-\gamma+\varsigma}(y)\log\eta_y dy  dxdt\rt|\\
 \lesssim &\bar Q^{\gamma+4}\mathfrak{E}_L(0) + \bar Q^{\gamma+2} \lt|\int_I \bar\rho v  \int_x^1 \bar\rho^{1-\gamma+\varsigma} \log \eta_y dy dx\big|_0^t \rt|\\
 &+  \bar Q^{\gamma+2}\lt|\int_0^t \int_I \bar\rho v  \int_x^1 \bar\rho^{1-\gamma+\varsigma} \eta_y^{-1} v_y dy dx dt\rt|\\
\le & \exp(C\bar Q^2)  \mathfrak{E}_L(0),
\end{split}
\eee
since
\bee\begin{split}
 & \lt|\int_I   \bar\rho v  \int_x^1 \bar\rho^{1-\gamma+\varsigma} \log \eta_y dy dx\big|_0^t \rt| \\
&\quad \lesssim   \bar Q^{-1} \lt|\int_I \bar\rho v  \lt(\int_x^1 \bar\rho^{1-\gamma+\varsigma}dy\rt)^{\frac 1 2} \lt(\int_x^1  \bar\rho^{1-\gamma+\varsigma}( \eta_y -1)^2 dy\rt)^{\frac 1 2 } dx(t) \rt|\\
&\qquad +¡¡\bar Q^{-1}\lt|\int_I \bar\rho v_0 \lt(\int_x^1 \bar\rho^{1-\gamma+\varsigma}dy\rt)^{\frac 1 2} \lt(\int_x^1  \bar\rho^{1-\gamma+\varsigma}( \eta_{0y} -1)^2 dy\rt)^{\frac 1 2 } dx  \rt|\\
 &\quad \le  \delta\int_I \bar{\rho}^{1-\gamma+\varsigma}(\eta_x -1)^2dx(t) + {C(\delta)}{\bar Q^2} \int_I \bar\rho^{2+\varsigma}  v^2 dx (t) + \exp(C\bar Q^2)  \mathfrak{E}_L(0),
\end{split}
\eee
and by using  Lemma \ref{lem33},
\bee\begin{split}
 &\lt|\int_0^t \int_I \bar\rho v  \int_x^1 \bar\rho^{1-\gamma+\varsigma} (\eta_y)^{-1} v_y dy dx dt\rt|\\
 \lesssim & \bar Q \lt|\int_0^t \int_I\bar\rho v \lt(\int_x^1 \bar\rho^{2-2\gamma+2\varsigma}dy\rt)^{\frac 1 2} \lt(\int_x^1    v_y^2 dy \rt)^{\frac 1 2}  dxdt \rt|\\
 \lesssim & \bar Q\int_0^t \int_I   v_x^2 dxdt + \bar Q \int_0^t\int_I\bar\rho^{3-\gamma+2\varsigma}
  v^2  dxdt \\
 \lesssim & \bar Q\int_0^t \int_I   v_x^2 dxdt.
\end{split}
\eee
Thus \ef{30B} is shown. To proceed on the energy estimates, we need to estimate $|\eta_x^{-1}v_x|$. We integrate \ef{8} in $(x,1)$ to get
 \be\label{27}\begin{split}
 \lt|\frac{v_x}{\eta_x} \rt| = &\lt|\frac{\bar\rho^{\gamma}}{\e^2}(\eta_x^{-\gamma}-1) - \int_x^1\bar\rho v_t dy\rt| \\
 \lesssim &\frac{\bar\rho^{\gamma}}{\e^2}\bar Q^{\gamma+1}|\eta_x-1| + \bar\rho^{\frac \gamma 2} \lt(\int_I \bar\rho v_t^2 dy\rt)^{\frac 1 2},
 \end{split}
 \ee
 due to \ef{oo} and \ef{11}.
Then, we multiply \ef{8} by $\e^2 v_t$ and integrate over $I$ to get
\be\label{vtl2l2}
\begin{split}
\e^2\int_I & \bar\rho  v_t
^2 dx + \frac{d}{dt}\int_I\lt[\bar\rho^\gamma\lt(\frac{1}{\gamma-1}(\eta_x)^{1-\gamma}
+\eta_x-\frac{\gamma}{\gamma-1}\rt)+\frac{\e^2}{2}  \eta_x^{-1}v_x^2\rt]dx\\
 &  =\frac{\e^2}{2}
\int_I \eta_x^{-2}v_x^3 dx \lesssim \e^2 \lt\|\frac{v_x}{\eta_x}\rt\|_{L^\infty} \int_I \eta_x^{-1}v_x^2 dx \\
&\le \frac{\e^2} 2 \int_I   \bar\rho v_t^2 dx + C\lt(\bar Q^{\gamma+2}+ \e^2 \int_I\eta_x^{-1} v_x^2 dx\rt)\int_I\eta_x^{-1} v_x^2 dx,
\end{split}
\ee
where \ef{27} and \ef{11} were used in the last inequality.
Thanks to the Gronwall inequality, we derive
\be\label{C18pp}
\begin{split}
  \e^2\int_0^t\int_I &   \bar\rho v_t^2 dxdt + \int_I [\bar\rho^\gamma(\eta_x-1)^2 + \e^2v_x^2 ]dx(t)\\
  \lesssim & \exp \lt(\int_0^t \int_I\eta_x^{-1} v_x^2 dx dt\rt) \lt(\mathfrak{E}_L(0)+\bar Q^{\gamma+2}   \int_0^t \int_I\eta_x^{-1} v_x^2 dx dt\rt)\\
    \lesssim & \exp(C\bar Q^2)  \mathfrak{E}_L(0).
\end{split}
\ee
Note that, with the aid of \ef{30B}, we have
\bee
\begin{split}
I_1:= & \int_0^t   (1+  s)^{\lambda-1}
\int_I (\eta_x -1)^2dxds\\
 \lesssim  & \int_0^t (1+s)^{-1}  \lt((1+s)^{\lambda}\int_I \bar\rho^{1-\gamma+\varsigma}(\eta_x  -1)^2 dx\rt)^{\frac{\ga }{2\ga-1-\varsigma}}  \\
&\qquad\times\lt((1+s)^{\lambda}\int_I \bar\rho^{\gamma}(\eta_x -1)^2dx\rt)^{\frac{\gamma-1-\varsigma}{2\ga-1-\varsigma}}  ds
  \\
  \le  & \da \int_0^t (1+s)^{\lambda}\int_I\bar\rho^{\gamma}(\eta_x -1)^2dx ds\\
   &+ \frac C \da \int_0^t (1+s)^{\lambda-\frac{2\gamma-1-\varsigma}{\gamma} }\int_I \bar\rho^{1-\gamma+\varsigma}(\eta_x  -1)^2 dxds\\
  \le  & \frac{\da}{\e^2} \int_0^t (1+s)^{\lambda}\int_I\bar\rho^{\gamma}(\eta_x -1)^2dx ds + \frac C \da \exp(C\bar Q^2)  \mathfrak{E}_L(0),
 \end{split}
\eee
  Taking $\beta=0$ in \ef{20} and multiplying the resulting equality by $(1+t)^\lambda$, $\lambda\in (0,\frac{\gamma-1-\varsigma}{\gamma})$, then integrating over $I$, we obtain
\be\label{33}
\begin{split}
   (1+t)^{\lambda}& \int_I (\eta_x -1)^2  dx(t) +\frac{1}{\e^2} \int_0^t (1+s)^{\lambda}
\int \bar\rho^{\gamma }(\eta_x -1)^2dxds \\
 \lesssim & \exp(C\bar Q^2)\lt(\mathfrak{E}_L(0) +   \int_0^t (1+s)^{\lambda} \int_I    v_x^2 dx ds+ (1+t)^\lambda\int_I \bar\rho v^2 dx(t)\rt.\\
 &\lt. +  \int_0^t (1+s)^{\lambda-1}
\int_I (\eta_x -1)^2dxds\rt)\\
\lesssim &  \exp(C\bar Q^2)  \mathfrak{E}_L(0),
 \end{split}
\ee
 where  Lemma  \ref{lem3} and the estimate for $I_1$ are used. Next, with \ef{21} and \ef{33}, we integrate $(1+t)^{1+\lambda}\cdot\ef{9}$, $\lambda\in (0,\frac{\gamma-1-\varsigma}{\gamma})$ in $I\times(0,t)$ to get
\be\label{34}
\begin{split}
    (1+t)^{1+\lambda} & \int_I\lt[ \bar\rho v^2 + \frac{\bar\rho^\gamma}{\e^2} (\eta_x-1)^2\rt] dx(t) +\int_0^t (1+s)^{1+\lambda}\int_I  v_x^2 dx   ds \\
   \lesssim & \exp(C\bar Q^2)\lt[  \mathfrak{E}_L(0) +\int_0^t(1+s)^\lambda\int_I\lt( \bar\rho v^2 +\frac{ \bar\rho^\gamma}{\e^2} (\eta_x-1)^2\rt) dxdt\rt]\\
  \lesssim & \exp(C\bar Q^2)  \mathfrak{E}_L(0) .
 \end{split}
\ee
Similarly as \ef{33}, we take $\beta=-\tilde{\beta}$, $\tilde{\beta}> 1$, in \ef{20}, and multiply the resulting equality by $(1+t)^{1+\lambda}$, $\lambda\in [0,\frac{\gamma-1-\varsigma}{\gamma})$, and integrate to get
\be\label{decay3}
\begin{split}
   (1+t)^{1+\lambda}& \int_I \bar\rho^{\tilde\beta}(\eta_x -1)^2  dx(t) +\frac{1}{\e^2} \int_0^t (1+s)^{1+\lambda}
\int \bar\rho^{\gamma+\tilde\beta }(\eta_x -1)^2dxds \\
 \lesssim & \exp(C\bar Q^2)\lt(\mathfrak{E}_L(0) +   \int_0^t (1+s)^{1+\lambda} \int_I    v_x^2 dx ds+ (1+t)^{1+\lambda}\int_I \bar\rho v^2 dx(t)\rt.\\
 &\lt. +  \int_0^t (1+s)^{\lambda}
\int_I \bar\rho^{\tilde\beta}(\eta_x -1)^2dxds\rt)\\
\lesssim &  \exp(C\bar Q^2)  \mathfrak{E}_L(0),
 \end{split}
\ee
where the last inequality holds for $\lambda=0$ and $\tilde\beta=1+\varsigma$ by using  \ef{30B} and \ef{34}, moreover, for $\lambda=\frac{\gamma-1}{\gamma}-\varsigma$ and $\tilde\beta=\gamma$ by applying \ef{33} and \ef{34}. Then we integrate $(1+t)^{1+\lambda}$ times  the first equality of     $\ef{vtl2l2}$, $\lambda\in (0,\frac{\gamma-1-\varsigma}{\gamma})$ and use the Gronwall inequality  to get
\be\label{decay4}
\begin{split}
\e^2\int_0^t & (1+s)^{1+\lambda}\int_I   \bar\rho  v_s
^2 dxds + (1+t)^{1+\lambda}\int_I\lt[\bar\rho^\gamma(\eta_x-1)^2+ \e^2   v_x^2\rt]dx (t)\\
 &  \lesssim \exp\lt\{\int_0^t\int_I \eta_x^{-1} v_x^2\rt\} \exp(C\bar Q^2)\lt(\mathfrak{E}_L(0) +   \int_0^t (1+s)^{1+\lambda} \int_I    v_x^2 dx ds  \rt.\\
 &\lt. +  \int_0^t (1+s)^{\lambda}
\int_I \bar\rho^{\gamma}(\eta_x -1)^2dxds\rt)\\
&\lesssim \exp(C\bar Q^2) \mathfrak{E}_L(0),
\end{split}
\ee
 by employing \ef{11}, \ef{33} and \ef{34}.
 Thus we prove \ef{30} by combining  \ef{C18pp}, \ef{33}, \ef{34}, \ef{decay4}  and using Lemma \ref{lem33}.

Next, for $\gamma\ge 2$, we have
\bee
\int_I \bar\rho^{2\gamma-2}(\eta_x-1)^2\lesssim \int_I \bar\rho^{\gamma}(\eta_x-1)^2 dx\le (1+t)^{-\frac{2\gamma-1}{\gamma}+\varsigma}\exp(C\bar Q^2) \mathfrak{E}_L(0)
\eee
by \ef{34}, while for $\gamma<2$, one gets
\bee
\begin{split}
\int_I \bar\rho^{2\gamma-2}(\eta_x-1)^2&\lesssim \lt(\int_I \bar\rho^{1-\gamma+\varsigma}(\eta_x-1)^2 dx\rt)^{\frac{2-\gamma}{2\gamma-1-\varsigma}}\lt(\int_I \bar\rho^{\gamma}(\eta_x-1)^2 dx\rt)^{\frac{3\gamma-3-\varsigma}{2\gamma-1-\varsigma}}\\
&\le (1+t)^{-\frac{3\gamma-3}{\gamma}+\varsigma}\exp(C\bar Q^2)\mathfrak{E}_L(0),
\end{split}
\eee
by \ef{30B} and \ef{34}.
Thus we derive from Lemma \ref{lem33} that
\be\label{38}
\begin{split}
\int_I(\eta-x)^2 dx&\lesssim \int_I \bar\rho^{2\gamma-2}(\eta_x-1)^2 dx\\
&\le C(1+t)^{-\frac{\min(2\gamma-2,\gamma)+\gamma-1}{\gamma}+\varsigma}\exp(C\bar Q^2) \mathfrak{E}_L(0).
\end{split}
\ee
It follows from \ef{33} and \ef{38} that
\be\label{39}
\begin{split}
\|\eta-x\|_{L^\infty}^2&=\lt\|\int_0^x[(\eta(y)-y)^2]_ydy\rt\|_{L^\infty}^2\\
&\lesssim \|\eta-x\|_{L^2}\|\eta_x-1\|_{L^2}\\
&\lesssim (1+t)^{-\frac{\min( \gamma-1,\gamma/2)+\gamma-1}{\gamma}+\varsigma}\exp(C\bar Q^2) \mathfrak{E}_L(0).
\end{split}
\ee
Therefore, by   \ef{38} and \ef{39}, we obtain \ef{36}.
\hfill$\Box$

\begin{lem}\label{lem4}  For  $t\ge 0$, we have
\be
\label{24}
\begin{split}
   (1& +t)^{\frac{2\gamma-1}{\gamma}-\varsigma}\e^4\int_I   \bar\rho v_t^2   dx (t) + \e^4 \int_0^t(1+s)^{\frac{2\gamma-1}{\gamma}-\varsigma}\int_I (v_s^2 +  {v_{sx}^2})dxds \\
   &  \lesssim  \exp(C\bar Q^2) (\mathfrak{E}_L(0)+ \mathfrak{E}_L^2(0)) ,
 \end{split}
 \ee
where $\varsigma \in (0,\frac{ \gamma-1}{\gamma})$ is an arbitrary constant.
\end{lem}
{\bf Proof.}
 We apply $\partial_t$ to \ef{8} to get
\be\label{25}
\bar\rho v_{tt}-\frac\gamma{\e^2} \lt( {\bar\rho^\gamma}{\eta_x^{-\gamma-1}}v_x\rt)_x= \lt({\eta_x^{-1 }}{v_{tx}} - {\eta_x^{-2 }}{v_x^2}\rt)_x.
\ee
Integrating    $\e^4v_t\cdot\ef{25}$ over $I$ gives
\be\label{26} \begin{split}
  \frac{\e^4}{2}\frac{d}{dt} & \int_I   \bar\rho v_t^2    dx+   \e^4\int_I {\eta_x^{-1 }}   {v_{tx}^2}dx \\
  =& -\gamma\e^2\int_I \bar\rho^\gamma {\eta_x^{-\gamma-1}}v_x v_{tx} dx + \e^4\int_I   {\eta_x^{-2 }} {v_x^2}v_{tx}dx \\
  \le  & \frac {\e^4} 2 \int_I {\eta_x^{-1 }}   {v_{tx}^2}dx+ C\lt(\|\eta_x^{-1}\|_{L^\infty_{x,t}}^{2\gamma}+ \e^4\| \eta_x^{-1}{v_x}\|_{L^\infty}^2  \rt) \int_I  \eta_x^{-1} {v_x^2} dx.
\end{split}\ee
Submit  \ef{27} into \ef{26}, then use \ef{11} to yield
\be\label{28}\begin{split}
  \e^4\frac{d}{dt}\int_I    \bar\rho v_t^2   dx+  \e^4 \int_I   \eta_x^{-1}   {v_{tx}^2}dx
   \lesssim      \lt[ \bar Q^{2\gamma+4}+ \e^4\int_I \bar\rho v_t^2 dx   \rt]  \int_I   \eta_x^{-1} {v_x^2} dx    .
\end{split}\ee
By the Gronwall inequality , \ef{7A} and Lemma \ref{lem2}, one gets, for $t\ge 0$,
\be \label{29}\begin{split}
\e^4 \int_I  & \bar\rho v_t^2   dx (t) +  \e^4\int_0^t\int_I   \eta_x^{-1} {v_{tx}^2}dxdt\\
 &\lesssim  \exp\lt( C \int_0^t\int_I  \eta_x^{-1}v_x^2 dx dt\rt) \lt(\e^4\int_I \bar\rho v_t^2(x,0)dx +  \bar Q^{2\gamma+2}\int_0^t \int_I  \eta_x^{-1}v_x^2 dx dt\rt)\\
  &\lesssim \exp(C\bar Q^2)(1+\bar Q^{2\gamma+4})\mathfrak{E}_L(0) \lesssim \exp(C\bar Q^2)   \mathfrak{E}_L(0).
\end{split}\ee
This estimate, with the help of \ef{27} and \ef{11} ,  implies
\bee
\e^2|v_x(x,t)|  \lesssim  \bar Q^{\gamma+3} +  \bar Q^{\gamma+3} \exp(C\bar Q^2)\mathfrak{E}_L(0)\lesssim \exp(C\bar Q^2)\mathfrak{E}_L(0) .
\eee
   for any  $(x,t) \in I \times [0,\iy)$.

In order to show the decay rates, we multiply \ef{28} by $1+t$ and integrate with respect to the time variable to get, for $t\ge 0$,
\be\label{29A}
\begin{split}
  (1 +t)\e^4 &\int_I   \bar\rho v_t^2  dx (t)+  \e^4\int_0^t(1+s)\int_I   \eta_x^{-1} {v_{sx}^2}dxds\\
   \lesssim & \e^4\int_I \bar\rho v_t^2(x,0) dx + \e^4\int_0^t\int_I \bar\rho v_t^2  dxdt\\
   & +
  \lt[ \bar Q^{2\gamma+4}+\e^4 \sup_{s\in[0,t]}\int_I \bar\rho v_s^2 dx(s)   \rt]  \int_0^t (1+s)\int_I \eta_x^{-1}  {v_x^2} dx ds\\
  \lesssim & \mathfrak{E}(0) + \e^4\int_0^t\int_I   {v_{tx}^2}dxdt +     \exp(C\bar Q^2)  \mathfrak{E}_L^2(0)\\
  \lesssim & \exp(C\bar Q^2) (\mathfrak{E}_L(0)+ \mathfrak{E}_L^2(0)  ) ,
\end{split}
\ee
with the help of \ef{21}, \ef{29}  and Lemma \ref{lem33}.

In a similar way, we multiply \ef{28} by $(1+t)^{1+\lambda} $,  $\lambda\in (0,\frac{\gamma-1-\varsigma}{\gamma})$ with $\varsigma \in (0,\gamma-1)$, and integrate to get
\be
\label{35}
\begin{split}
    (1 +t)^{1+\lambda} &\e^4\int_I  \bar\rho v_t^2   dx(t) +  \e^4\int_0^t(1+s)^{1+\lambda}\int_I    {v_{sx}^2} dxds \\
   \lesssim  &  \bar Q\lt(\e^4\int_I  \bar\rho v_t^2(x,0)   dx +\e^4 \int_0^t(1+s)^{ \lambda}\int_I \bar\rho   {v_{s }^2} dxds \rt)\\
   &+ \bar Q^2\e^4\int_0^t(1+s)^{1+\lambda} \lt[ \bar Q^{2\gamma+4}+ \int_I \bar\rho v_t^2 dx   \rt]  \int_I   {v_x^2} dxds\\
   \lesssim & \bar Q\mathfrak{E}(0) + \bar Q\e^4\int_0^t(1+s) \int_I     {v_{sx }^2} dxds \\
   &+     \exp(C\bar Q^2)\mathfrak{E}(0) \int_0^t(1+s)^{1+\lambda}  \int_I   {v_x^2} dxds\\
   \lesssim & \exp(C\bar Q^2) (\mathfrak{E}_L(0)+ \mathfrak{E}_L^3(0)  ),
 \end{split}
 \ee
where \ef{11}, \ef{29A}, \ef{30} and Lemma \ref{lem33} are used. \hfill$\Box$

Recalling $(\bar\rho^\gamma)_x=-\bar\rho$, we rewrite \ef{19} as follows:
\be\label{40}
 (\eta_x^{-1}\eta_{xx})_t+\frac{\gamma}{\e^2}\bar{\rho}^{\gamma}\eta_x^{-\gamma-1}\eta_{xx}+ \frac{\bar{\rho}}{\e^2}(\eta_x^{-\gamma}-1)-\bar{\rho}v_t=0.
\ee

\begin{lem}\label{lem7} For any $t\ge 0$, we have
\be\label{41}
 \int_I\bar{\rho}^{\gamma-1+\varsigma}\eta_{xx}^2dx(t)
 +\frac{1}{\e^2}\int_0^t\int_I\bar{\rho}^{2\gamma-1+\varsigma}\eta_{xx}^2dxds\lesssim \exp(C\bar Q^2)  \mathfrak{E}_L(0),
\ee
\be\label{42}
\begin{split}
 (1+t)&\int_I\bar{\rho}^{2\gamma-1+\varsigma}\eta_{xx}^2dx(t)+\frac{1}{\e^2}\int_0^t(1+s)\int_I\bar{\rho}^{3\gamma-1+\varsigma}\eta_{xx}^2dxds\\
 & \lesssim \exp(C\bar Q^2)  (\mathfrak{E}_L(0)+ \mathfrak{E}_L^3(0)),
 \end{split}
\ee
\be\label{44}
(1+t)^{\frac{2\gamma-1}{\gamma}-\varsigma}\e^4\|v_x(\cdot,t)\|_{L^{\infty}}^2\lesssim  \exp(C\bar Q^2)  (\mathfrak{E}_L(0)+ \mathfrak{E}_L^3(0)),
\ee
where $\varsigma\in(0,\ga-1)$ is an arbitrary constant.
\end{lem}
{\bf Proof.} Multiplying \ef{40} by $\bar{\rho}^{\gamma-1+\varsigma}\eta_x^{-1}\eta_{xx}$ and integrating the resulting equality in $I$, we obtain
\bee\begin{split}
 \frac{1}{2}\frac{d}{dt}&\int_I\bar{\rho}^{\gamma-1+\varsigma}\eta_x^{-2}\eta_{xx}^2dx
 +\frac{\gamma}{\e^2}\int_I\bar{\rho}^{2\gamma-1+\varsigma}\eta_x^{-\gamma-2}\eta_{xx}^2dx\\
= & \int_I\frac{\bar{\rho}^{\gamma+\varsigma}}{\e}\eta_{x}^{-1}\eta_{xx}\lt[ \frac{1}{\e}(\eta_x^{-\gamma}-1)+\e v_t\rt]dx \\
\le & \frac{\delta}{\e^2}\int_I\bar{\rho}^{2\gamma-1+\varsigma}\eta_{xx}^2dx+ C(\delta)\bar Q^{\gamma+2}\int_I\bar{\rho}^{1+\varsigma}\lt[\frac{1}{\e^2}(\eta_x-1)^2+\e^2 v_t^2\rt]dx.
\end{split}
\eee
Choosing $\delta$ to be small and integrating with respect to $t$, then using \ef{30B} and \ef{30}, we get \ef{41}.

Next, we integrate $(1+t)\bar{\rho}^{2\gamma-1+\varsigma}\eta_x^{-1}\eta_{xx}\cdot \ef{40}$ over $I\times(0,t)$ to achieve
\bee\begin{split}
  (1+t)&\int_I\bar{\rho}^{2\gamma-1+\varsigma}\eta_x^{-2}\eta_{xx}^2dx(t)
+\frac{1}{\e^2}\int_0^t(1+s)\int_I\bar{\rho}^{3\gamma-1+\varsigma}\eta_x^{-\gamma-2}\eta_{xx}^2dxds\\
\lesssim & \int_I \bar{\rho}^{\gamma-1+\varsigma}\eta_{0x}^{-2}\eta_{0xx}^2dx+ \int_0^t\int_I\bar{\rho}^{2\gamma-1+\varsigma}\eta_x^{-2}\eta_{xx}^2dxdt\\
&+\bar Q^{\gamma+2}\int_0^t(1+s)\int_I\bar{\rho}^{2\gamma+\varsigma}\frac{|\eta_{xx}|}{\e}\lt[\frac{|\eta_x-1|}{\e}+\e |v_s|\rt]dxds\\
\lesssim & \exp(C\bar Q^2)  \mathfrak{E}(0)+\frac{\delta}{\e^2}\int_0^t(1+s)\int_I\bar{\rho}^{3\gamma-1+\varsigma}\eta_x^{-\gamma-2}\eta_{xx}^2dxds\\
&+C(\delta)\exp(C\bar Q^2) \int_0^t(1+s)\int_I\bar{\rho}^{\gamma+1+\varsigma}\lt[\frac{1}{\e^2}(\eta_x-1)^2+\e^2v_s^2\rt]dxds.
\end{split}
\eee
Absorbing the small term with $\delta$ on the right hand side,  using \ef{30} and \ef{decay3}, we achieve \ef{42}.
Similarly, we multiply \ef{40} by $(1+t)^{1+\sigma}\bar\rho^{3\gamma-2}\eta_x^{-1}\eta_{xx}$, $\sigma>0$ is a constant to be determined, and integrate with respect to the spatial variable to get
\bee\begin{split}
  (1+t)^{1+\sigma}&\int_I\bar{\rho}^{3\gamma-2}\eta_x^{-2}\eta_{xx}^2dx(t)
+\frac{1}{\e^2}\int_0^t(1+s)^{1+\sigma}\int_I\bar{\rho}^{4\gamma-2}\eta_x^{-\gamma-2}\eta_{xx}^2dxds\\
\lesssim & \int_I \bar{\rho}^{3\gamma-2}\eta_{0x}^{-2}\eta_{0xx}^2dx+ \int_0^t(1+s)^{ \sigma}\int_I\bar{\rho}^{3\gamma-2}\eta_x^{-2}\eta_{xx}^2dxds\\
&+\bar Q^{2\gamma+4}\int_0^t(1+s)^{1+\sigma}\int_I\bar{\rho}^{2\gamma} \lt[\frac{1}{\e^2}(\eta_x-1)^2+ \e^2v_s^2\rt]dxds,
\end{split}
\eee
where the last integral on the right hand side is estimated by \ef{30} and \ef{decay3}, and the second integral on the right hand side is estimated by using \ef{41} as follows:
\bee\begin{split}
        \int_0^t &  (1+s)^{ \sigma}\int_I\bar{\rho}^{3\gamma-2}\eta_x^{-2}\eta_{xx}^2dxds\\
       \lesssim &    \exp\{C\bar Q^2\} \int_0^t\lt[(1+s)^{1+\sigma}\int_I \bar\rho^{4\gamma-2}\eta_x^{-\gamma-2}\eta_{xx}^2 dx\rt]^\theta  \\
       &\qquad \times(1+s)^{\sigma(1-\theta)-\theta}\lt(\int_I \bar\rho^{\gamma-1+\varsigma}\eta_x^{-2}\eta_{xx}^2 dx\rt)^{1-\theta}ds\\
       \le &  \delta \int_0^t (1+s)^{1+\sigma}\int_I \bar\rho^{4\gamma-2}\eta_x^{-\gamma-2}\eta_{xx}^2 dxds \\
       &\quad+ C(\delta)\exp\{C\bar Q^2\} \int_0^t(1+s)^\frac{\sigma(1-\theta)-\theta}{1-\theta}ds\sup_{t\ge 0}\int_I \bar\rho^{\gamma-1+\varsigma}\eta_x^{ -2}\eta_{xx}^2 dx (t)
\end{split}
\eee
for $\theta=\frac{2\gamma-1-\varsigma}{3\gamma-1-\varsigma}\in(\frac 1 2,1)$ since $\varsigma\in (0,\gamma-1)$,  and $0<\sigma<\frac{2\theta-1}{1-\theta}=\frac{\gamma-1-\varsigma}{\gamma}$.
Thus, from \ef{27}, we find that
\be\label{vxlf}\begin{split}
 \e^4\lt\| v_x \rt\|_{L^\infty}^2 \lesssim & \bar Q^{2(\gamma+1)}\lt\|\bar\rho^{2\gamma}(\eta_x -1)^2\rt\|_{L^\infty}+  \bar Q^2\e^4\int_I \bar\rho v_t^2 dx\equiv J_1+ J_2,
 \end{split}
 \ee
 where
 \bee\begin{split}
 J_1&\lesssim   \bar Q^{2(\gamma+1)}\lt|\int_x^1[\bar\rho^\gamma(\eta_y-1)][\bar\rho^\gamma(\eta_y-1)]_y dy\rt|\\
 &\lesssim \bar Q^{2(\gamma+1)}\int_I \bar\rho^\gamma|\eta_x-1|(\bar\rho|\eta_x-1|+\bar\rho^{\gamma-1}|\eta_{xx}|) dx\\
 &\lesssim \bar Q^{2(\gamma+1)}\lt(\int_I \bar\rho^{\gamma}(\eta_x-1)^2dx + \int_I \bar\rho^{3\gamma-2} \eta_{xx}^2dx\rt)\\
 &\lesssim (1+t)^{-\frac{2\gamma-1}{\gamma}+\varsigma}\exp(C\bar Q^2)  (\mathfrak{E}(0)+ \mathfrak{E}^3(0)) .
\end{split}
 \eee
Therefore, we show \ef{44} by \ef{24},  \ef{vxlf} and the estimate for $J_1$.

\hfill$\Box$

\begin{lem}\label{lem8}  For any $t\ge0$, we have
\be\label{49}
\e^4\|v_{xx}(\cdot,t)\|_{L^2}^2(t)\lesssim (1+t)^{-\frac{1}{\gamma} \min{ (1,\gamma-1   ) }+\varsigma}\exp(C\bar Q^2)  (\mathfrak{E}_L(0)+ \mathfrak{E}_L^3(0)),
\ee
\be\label{49-1}
\e^2\|\eta_{xx}(\cdot,t)\|_{L^2}^2(t)\lesssim (1+t)^{\frac{\gamma-1}{\gamma}+\varsigma}\exp(C\bar Q^2)  (\mathfrak{E}_L(0)+ \mathfrak{E}_L^3(0)),
\ee
where $\varsigma\in(0,\frac{1}{\gamma}\min (\gamma-1,1))$ is an arbitrary constant.
\end{lem}
{\bf Proof.} To estimate $v_{xx}$, we integrate \ef{8} from $x$ to $1$ and multiply the resulting equation by $\e^2\eta_x$ and then differentiate with respect to $x$ to get
\be\label{50}\begin{split}
 \e^2 v_{xx}= & \lt[\bar{\rho}^{\gamma}(\eta_x^{-\gamma}-1)-\e^2\int_x^1\bar{\rho}v_tdy\rt]\eta_{xx}\\
& + \eta_x\lt[-\gamma\bar{\rho}^{\gamma}\eta_x^{-\gamma-1}\eta_{xx}- \bar{\rho}(\eta_x^{-\gamma}-1)+\e^2\bar{\rho}v_t\rt],
\end{split}\ee
which gives
\be\label{51}
\e^4\|v_{xx}\|_{L^2}^2\lesssim \int_I\bar{\rho}^{\gamma}[(\eta_x-1)^2+\eta_{xx}^2]dx+\e^4\int_I\bar{\rho}v_t^2dx+\int_I\bar{\rho}^2(\eta_x-1)^2dx.
\ee
Note that, for $\gamma\le 2$ and $\varsigma\in (0,\frac{\gamma-1}{\gamma})$,
\bee\begin{split}
\int_I\bar{\rho}^2(\eta_x-1)^2dx&\lesssim \int_I\bar{\rho}^{\gamma}(\eta_x-1)^2dx\\
&\lesssim (1+t)^{-\frac{2\gamma-1}{\gamma}+\varsigma}\exp(C\bar Q^2)  \mathfrak{E}_L(0) ,
\end{split}\eee
and for $\gamma > 2$ and $\varsigma\in (0,\frac{\gamma-1}{\gamma})$,
\bee\begin{split}
\int_I\bar{\rho}^2(\eta_x-1)^2dx & \lesssim \lt(\int_I\bar{\rho}^{1-\gamma+\varsigma}(\eta_x-1)^2dx\rt)^{\frac{\gamma-2}{2\gamma-1-\varsigma}}
\lt(\int_I\bar{\rho}^\gamma(\eta_x-1)^2dx\rt)^{\frac{\gamma+1-\varsigma}{2\gamma-1-\varsigma}}\\
& \lesssim (1+t)^{-\frac{\gamma+1}{\gamma}+\varsigma}\exp(C\bar Q^2)  \mathfrak{E}_L(0) ,
\end{split}\eee
in virtue of \ef{30B} and \ef{30}. Thus
\be\label{52}
\int_I\bar{\rho}^2(\eta_x-1)^2dx \lesssim (1+t)^{-1-\frac{\min(\gamma-1,1)}{\gamma}+\varsigma}\exp(C\bar Q^2)  \mathfrak{E}_L(0).
\ee
Next, we multiply \ef{40} by $(1+t)^{\sigma}\bar{\rho}^{\gamma}\eta_x^{-1}\eta_{xx}$, and integrate to get
\be\label{53}\begin{split}
 \frac{1}{2}(1+& t)^{\sigma}\int_I\bar{\rho}^{\gamma}\eta_x^{-2}\eta_{xx}^2dx(t)
 +\frac{\gamma}{\e^2}\int_0^t(1+s)^{\sigma}\int_I\bar{\rho}^{2\gamma}\eta_x^{-\gamma-2}\eta_{xx}^2dxds\\
\lesssim & \int_I \bar\rho^\gamma\eta_{0x}^{-2}\eta_{0xx}^2dx + \int_0^t(1+s)^{\sigma-1}\int_I\bar{\rho}^{\gamma}\eta_{xx}^2dxds\\
&+\int_0^t(1+s)^{\sigma}\int_I\bar{\rho}^2\lt[\frac{1}{\e^2}(\eta_x-1)^2
+\e^2v_s^2\rt]dxds.
\end{split}\ee
Let $\sigma\in(0,\frac{1}{\gamma}\min (\gamma-1,1)-\varsigma]$  with $\varsigma\in(0,\frac{1}{\gamma}\min (\gamma-1,1))$. Note that, for $\gamma\le2$,
\bee\begin{split}
  \frac{1}{\e^2}\int_0^t(1+s)^{\sigma}\int_I\bar{\rho}^2(\eta_x-1)^2dxds
&\lesssim  \frac{1}{\e^2} \int_0^t(1+s)^{\sigma}\int_I\bar{\rho}^\gamma(\eta_x-1)^2dxds\\
&\lesssim  \exp(C\bar Q^2)   \mathfrak{E}_L(0) ,
\end{split}\eee
and for $\gamma>2 $,
 \bee\begin{split}
  \frac{1}{\e^2}\int_0^t(1& +s)^{\sigma}\int_I\bar{\rho}^2(\eta_x-1)^2dxds\\
\lesssim & \frac{1}{\e^2}\int_0^t(1+s)^{\sigma-\frac{\gamma-1-\varsigma}{\gamma} \frac{1-\varsigma}{\gamma-1-\varsigma}}\lt[  \int_I\bar{\rho}^{1+\varsigma}(\eta_x-1)^2dx
\rt]^{\frac{\gamma-2}{\gamma-1-\varsigma}}\\
& \times\lt[(1+s)^{\frac{\gamma-1-\varsigma}{\gamma}}
\int_I\bar{\rho}^{\gamma}(\eta_x-1)^2dx\rt]^{\frac{1-\varsigma}{\gamma-1-\varsigma}}ds\\
\le & \frac{1}{\e^2} \int_0^t(1+s)^{\frac{\gamma-1-\varsigma}{\gamma}}\int_I\bar{\rho}^\gamma(\eta_x-1)^2dxds\\
&+\frac{1}{\e^2}\int_0^t(1+s)^{(\sigma- \frac{1-\varsigma}{\gamma} )\frac{\gamma-1-\varsigma}{\gamma-2}} \int_I\bar{\rho}^{1+\varsigma}(\eta_x-1)^2dx ds  \\
\lesssim &  \exp(C\bar Q^2)   \mathfrak{E}_L(0) ,
\end{split}
\eee
by observing that $\frac{\gamma-1-\varsigma}{\gamma-2}>0$ and $\sigma- \frac{1-\varsigma}{\gamma} \le 0$, and employing \ef{30B} and \ef{30}. It follows that
\be\label{54}
\frac{1}{\e^2}\int_0^t(1  +s)^{\sigma}\int_I\bar{\rho}^2(\eta_x-1)^2dxds \lesssim  \exp(C\bar Q^2)   \mathfrak{E}_L(0) ,\quad \sigma\in(0,\frac{1}{\gamma}\min (\gamma-1,1)-\varsigma).
\ee
Similarly, we use \ef{41} and the interpolation to show that
\be\label{55}\begin{split}
  \int_0^t(1&+s)^{\sigma-1}\int_I\bar{\rho}^{\gamma}\eta_{xx}^2dxds\\
\lesssim & \int_0^t(1+s)^{\sigma-1}\lt[\int_I\bar{\rho}^{\gamma-1+\varsigma}\eta_{xx}^2dx\rt]^{\frac{\gamma-1+\varsigma}{\gamma}}\lt[\int_I\bar{\rho}^{2\gamma-1+\varsigma}\eta_{xx}^2dx\rt]^{\frac{1-\varsigma}{\gamma}}ds\\
\lesssim &
\int_0^t\int_I\bar{\rho}^{2\gamma-1+\varsigma}\eta_{xx}^2dxds
+\int_0^t(1+s)^{(\sigma-1)\frac{\gamma}{\gamma-1+\varsigma}}ds \int_I\bar{\rho}^{\gamma-1+\varsigma}\eta_{xx}^2dx\\
\lesssim &  \exp(C\bar Q^2)  (\mathfrak{E}_L(0)+ \mathfrak{E}_L^3(0)) ,
\end{split}\ee
when $(\sigma-1)\frac{\gamma}{\gamma-1+\varsigma}<-1$, that is, $\sigma<\frac{1-\varsigma}{\gamma}$.

Inserting \ef{54} and \ef{55} into \ef{53}, and using \ef {30}, we obtain
\be\label{56}
(1+t)^{\frac{1}{\gamma} \min(1,\gamma-1)-\varsigma}\int_I\bar{\rho}^{\gamma}\eta_{xx}^2dx(t)\lesssim \exp(C\bar Q^2)  (\mathfrak{E}_L(0)+ \mathfrak{E}_L^3(0)) .
\ee
Therefore, we conclude \ef{49} from \ef{30}, \ef{24}, \ef{51}, \ef{52} and \ef{56}.

Finally, multiplying \ef{40} by $\e^2\eta_x^{-1}\eta_{xx}$ and integrating over $I\times(0,t)$, we find
\bee\begin{split}
 \frac{\e^2}{2} \int_I & \eta_x^{-2}\eta_{xx}^2dx(t)+ {\gamma} \int_0^t\int_I\bar{\rho}^\gamma\eta_x^{-\gamma-2}\eta_{xx}^2dxdt \\
=& \e^2\int_I\eta_{0x}^{-2}\eta_{0xx}^2dx+ \int_0^t\int_I\eta_x^{-1} {\eta_{xx}} \lt[ {\bar{\rho}} (\eta_x^{-\gamma}-1)+\e^2\bar{\rho}v_t\rt]dxdt\\
\le & \mathfrak{E}_L(0)+\frac{\gamma}{2} \int_0^t\int_I\bar{\rho}^{\gamma}\eta_x^{-\gamma-2}\eta_{xx}^2dxdt+C  \bar Q^{2\gamma+4}\int_0^t  \int_I \bar{\rho}^{2-\gamma}[(\eta_x-1)^2+\e^4 v_t^2]dxdt,
\end{split}\eee
where the last term is estimated as follows: for $\varsigma\in (0,1)$,
\bee\begin{split}
\int_0^t & \int_I \bar{\rho}^{2-\gamma} (\eta_x-1)^2dxdt\\
 & \lesssim  \int_0^t\lt(\int_I\bar\rho^{1-\gamma+\varsigma}(\eta_x-1)^2dx\rt)^{\frac{2\gamma-2}{2\gamma-1-\varsigma}}
 \lt(\int_I\bar{\rho}^{\gamma}(\eta_x-1)^2dx\rt)^{\frac{1-\varsigma}{2\gamma-1-\varsigma}}dt\\
& \lesssim     (1+t)^{\frac{\gamma-1}{\gamma}+\varsigma}\exp(C\bar Q^2)  (\mathfrak{E}_L(0)+ \mathfrak{E}_L^3(0)) ,
\end{split}\eee
by using \ef{30B} and \ef{30}, and
$$\e^4\int_0^t\int_I \bar{\rho}^{2-\gamma} v_t^2 dxdt\lesssim\e^4\int_0^t\int_I\bar{\rho}^\gamma v_{tx}^2dxdt\lesssim \exp(C\bar Q^2)  (\mathfrak{E}_L(0)+ \mathfrak{E}_L^2(0)) ,$$
by applying Lemma \ref{hardy2} and \ef{24}. Therefore \ef{49-1} is shown. \hfill$\Box$
\vskip 2mm
\noindent{\bf Proof of Theorem 2.2, part (i)}: { First, we show the regularity in \ef{lregu} and the global existence of strong solution for each $\e$. From \ef{44}, \ef{49}, the boundary condition $v(0,t)=0$ and Lemma \ref{hardy2}, we obtain that
\be\label{ee1}
\|v\|_{ L^\infty(0,T;H^2(I))}\le C(\epsilon).
\ee
Note that from \ef{3}, \ef{eta0def}, \ef{lag} and $\ef{6}_2$, we get
\bee
\eta(0,t)=\eta_0(0) + \int_0^t v(0,s)ds=0,
\eee
and thus
\bee
|\eta(x,t)|=\left|\eta(0,t)+\int_0^x \eta_y dy\right|\lesssim \|\eta_x\|_{L^\infty(0,T;L^\infty(I))}.
\eee
Then we conclude from \ef{11}, \ef{49-1} and the above inequality that
\be\label{ee2}
\|\eta\|_{ L^\infty(0,T;H^2(I))} + \|\eta^{-1}\|_{L^\infty(0,T;L^\infty(I))}\le C(T,\epsilon).
\ee
Moreover, from $v_t(0,t)=0$, Lemma \ref{hardy2} and \ef{24}, we derive that
\be\label{ee3}
\|\sqrt{\bar\rho}v_t\|_{L^\infty(0,T;L^2(I))} + \| v_t\|_{L^2(0,T;H^1(I))}\le C(\epsilon).
\ee
Due to \ef{ee1}, \ef{ee2} and \ef{ee3}, we find that  each term in $\ef{6}_1$, that is, $\bar\rho v_t$, $\e^{-2}(\bar\rho^\gamma\eta_x^{-\gamma})_x=\e^{-2}(\gamma\bar\rho^{\gamma-1}\bar\rho_x\eta_x^{-\gamma} -\gamma\bar\rho^\gamma\eta_x^{-\gamma-1}\eta_{xx})$, $-\e^{-2}\bar\rho g$, $(v_x\eta_x^{-1})_x=\eta_x^{-1}v_{xx}
-\eta_x^{-2}\eta_{xx}v_x$, is bounded in $L^\infty(0,T;L^2(I))$, and the boundary terms in $\ef{6}_2$, $v$ and $v_x$  are bounded in $L^\infty(0,T;L^2(\partial I))$ by the trace theorem, where the bounds  only depend on $T$ and $\e$. Therefore, we can derive the global existence and uniqueness of the strong solution of \ef{6}, by applying the local existence result in Theorem 2.1 and the global in time bounds.

Next, the uniform estimates with respect to $\e\in(0,1]$ and $t\ge 0$ are shown as follows. The bound \ef{e} for the low-order energy is obtained by using \ef{7A}, \ef{11}, \ef{30B}, \ef{30}, \ef{24}, \ef{41} and \ef{49-1}. The estimates \ef{e1} of the spatial-time integrals are derived from \ef{30B}, \ef{30}, \ef{41} and \ef{42}. The decay estimates \ef{decay1} are achieved by concluding \ef{30}, \ef{36}, \ef{24}, \ef{44} and \ef{49}.
}

\subsection{Uniform higher-order estimates and global classical solutions}
\begin{lem}\label{lem9} For any $t\ge 0$, we have
\be\label{57}
\begin{split}
(1&+t)^{\frac{2\gamma-1}{\gamma}-\varsigma}\e^8\int_I\bar{\rho}v_{tt}^2dx(t)
+\e^8\int_0^t(1+s)^{\frac{2\gamma-1}{\gamma}-\varsigma}\int_Iv_{ssx}^2dxds\\
&+(1+t)^{\frac{2\gamma-1}{\gamma}-\varsigma}\e^8\|v_{tx}\|_{L^\infty}^2(t)\lesssim \exp(C\bar Q^2)  (\mathfrak{E}(0)+ \mathfrak{E}^7(0)),
\end{split}
\ee
where $\varsigma\in (0,\frac{\gamma-1}{\gamma})$ is an arbitrary   constant.
\end{lem}
{\bf Proof.} Differentiating \ef{25} in $t$, one gets
\be\label{58}
\bar{\rho}v_{ttt}-\frac{\gamma}{\e^2}(\bar{\rho}^{\gamma}\eta_x^{-\gamma-1}v_{tx}
- (\gamma+1) \bar{\rho}^{\gamma}\eta_x^{-\gamma-2}v_x^2)_x
= (\eta_x^{-1}v_{ttx}-3\eta_x^{-2}v_xv_{tx}+2\eta_x^{-3}v_x^3)_x
\ee
Multiplying this equality by $\e^8v_{tt}$ and integrate over $I$, we obtain
\be \label{59}\begin{split}
  \frac{\e^8}{2}\frac{d}{dt}&\int_I\bar{\rho}v_{tt}^2dx+ \e^8\int_I\eta_x^{-1}v_{ttx}^2dx\\
=& \int_I \e^4v_{ttx}\lt[ \e^4 (3\eta_x^{-2}v_xv_{tx}-2\eta_x^{-3}v_x^3)+\gamma\bar{\rho}^{\gamma}\e^2[(\gamma+1)\eta_x^{-\gamma-2}v_x^2-\eta_x^{-\gamma-1}v_{tx}]\rt]dx\\
\le & \delta \e^8\int_I \eta_x^{-1}v_{ttx}^2dx+C(\delta)\bar Q^{2\gamma+3} [(1+\e^4\|v_x\|_{L^{\infty}}^2)\e^4\int_Iv_{tx}^2dx
+(\e^4\|v_x\|_{L^\infty}^2+\e^8\|v_x\|_{L^\infty}^4)\int_Iv_x^2dx ]
\end{split}\ee
which gives, for $t\ge 0$,
\be\label{60}
\e^8\int_I\bar{\rho}v_{tt}^2dx(t)+ \e^8\int_0^t\int_Iv_{ttx}^2dxdt\lesssim  \exp(C\bar Q^2)  (\mathfrak{E}(0)+ \mathfrak{E}^7(0)),
\ee
by applying \ef{11}, \ef{21}, \ef{24} and \ef{44}.
Next, we integrate $(1+t)\cdot \ef{59}$ in $(0,t)$ to get
\be\label{61}\begin{split}
  (1+t& )\e^8\int_I\bar{\rho}v_{tt}^2dx(t)+ \e^8\int_0^t(1+s)\int_Iv_{ssx}^2dxds\\
\lesssim & \mathfrak{E}(0)+\e^8\int_0^t\int_I\bar{\rho}v_{tt}^2dxdt+\exp\{C\bar Q^2\}\lt[(1+ \e^4\|v_x\|^2_{L^\infty_{x,t}})\int_0^t(1+s)\int_I \e^4 v_{sx}^2 dxds\rt.\\
&\quad \lt.+  (1+ \e^8\|v_x\|^4_{L^\infty_{x,t}})\int_0^t(1+s)\int_I v_x^2 dxds\rt] \\
\lesssim & \exp(C\bar Q^2)  (\mathfrak{E}(0)+ \mathfrak{E}^7(0))+\e^8\int_0^t\int_Iv_{ttx}^2dxdt\\
\lesssim & \exp(C\bar Q^2)  (\mathfrak{E}(0)+ \mathfrak{E}^7(0)),
\end{split}\ee
by using \ef{11}, \ef{21},\ef{24}, \ef{44}, \ef{60} and Lemma \ref{lem33}.

To gain faster decay for $\int_I\bar{\rho}v_{tt}^2dx(t)$, we multiply the first equality in \ef{59} by $(1+t)^{\frac{2\gamma-1}{\gamma}-\varsigma}$  and integrate the resulting equality to get
\be\label{vtt}\begin{split}
  (1+t)^{\frac{2\gamma-1}{\gamma}-\varsigma} & \e^8\int_I\bar{\rho}v_{tt}^2dx(t)+\e^8\int_0^t(1+s)^{\frac{2\gamma-1}{\gamma}-\varsigma}\int_Iv_{ssx}^2dxds\\
\lesssim & \mathfrak{E}(0)+\e^8\int_0^t(1+s)^{\frac{\gamma-1}{\gamma}-\varsigma}\int_I\bar{\rho}v_{tt}^2dxds
\\
&\quad+\exp\{C\bar Q^2\}\lt[(1+ \e^4\|v_x\|^2_{L^\infty_{x,t}})\int_0^t(1+s)^{\frac{2\gamma-1}{\gamma}-\varsigma}\int_I \e^4 v_{sx}^2 dxds\rt.\\
&\quad \lt.+  (1+ \e^8\|v_x\|^4_{L^\infty_{x,t}})\int_0^t(1+s)^{\frac{2\gamma-1}{\gamma}-\varsigma}\int_I v_x^2 dxds\rt]\\
\lesssim & \exp(C\bar Q^2)  (\mathfrak{E}(0)+ \mathfrak{E}^7(0))+\e^8\int_0^t(1+s)\int_Iv_{ssx}^2dxds\\
\lesssim & \exp(C\bar Q^2)  (\mathfrak{E}(0)+ \mathfrak{E}^7(0)),
\end{split}
\ee
where \ef{11}, \ef{30}, \ef{24}, \ef{44}, \ef{61} and Lemma \ref{lem33} were employed.

Integrating \ef{8} from $x$ to $1$ and applying $\partial_t$ to the resulting equality, we get
\be\label{62}
 v_{tx}=v_x\lt[\frac{\bar{\rho}^{\gamma}}{\e^2}(\eta_x^{-\gamma}-1 -\gamma\eta_x^{-\gamma-1}) - \int_x^1\bar{\rho}v_tdy\rt]-\eta_x\int_x^1\bar{\rho}v_{tt}dy.
\ee
It is not hard to find
\be\label{63}
\begin{split}
 \e^8\|v_{tx}\|_{L^\infty}^2 & \le \e^4\|v_x\|_{L^\infty}^2[1+ \bar Q^{\gamma+1} + \e^4\int_I \bar\rho v_t^2 dx] +\bar Q^2\e^8\int_I\bar{\rho}v_{tt}^2dx\\
&\le (1+t)^{-\frac{2\gamma-1}{\gamma}+\varsigma}\exp(C\bar Q^2)  (\mathfrak{E}(0)+ \mathfrak{E}^7(0)),
\end{split}
\ee
by using \ef{11}, \ef{24}, \ef{44} and \ef{vtt}. \hfill$\Box$

\begin{lem}\label{lem10} For any $t\ge 0$, there holds
\be\label{64}
\e^2\int_I\bar{\rho}^{3\gamma-3+\varsigma}\eta_{xxx}^2dx(t)
+\int_0^t\int_I\bar{\rho}^{4\gamma-3+\varsigma}\eta_{xxx}^2dxdt\lesssim \exp(C\bar Q^2)  (\mathfrak{E}(0)+ \mathfrak{E}^7(0)),
\ee
\end{lem}
where $\varsigma\in (0,\frac{\gamma-1}{\gamma})$ is an arbitrary constant.

\noindent{\bf Proof.} Applying $\partial_x$ to \ef{40}, we get
\be\label{65}\begin{split}
 (\eta_x^{-1}\eta_{xx})_{xt} & +\frac{\gamma}{\e^2}\bar{\rho}^{\gamma}\eta_x^{-\gamma}(\eta_x^{-1}\eta_{xx})_x \\
&=-\frac{\gamma}{\e^2}(\bar{\rho}^{\gamma}\eta_x^{-\gamma})_x\eta_x^{-1}\eta_{xx}- \frac{1}{\e^2}[\bar{\rho}(\eta_x^{-\gamma}-1)]_x+(\bar{\rho}v_t)_x.
\end{split}\ee
Multiplying \ef{65} by $ \e^2(\eta_x^{-1}\eta_{xx})_x\bar{\rho}^{3\gamma-3+3\varsigma} $ and integrating over $I\times(0,t)$, we discover
\bee\begin{split}
LHS& :=  \e^2\int_I\bar{\rho}^{3\gamma-3+3\varsigma}|(\eta_x^{-1}\eta_{xx})_x|^2dx(t)+\gamma\int_0^t\int_I\bar{\rho}^{4\gamma-3+3\varsigma}\eta_x^{-\gamma}|(\eta_x^{-1}\eta_{xx})_x|^2dxdt\\
&\lesssim \e^2\int_I\bar{\rho}^{3\gamma-3+3\varsigma}|(\eta_{0x}^{-1}\eta_{0xx})_x|^2dx\\
&\ \ + \int_0^t\int_I\bar{\rho}^{3\gamma-3+3\varsigma}|(\eta_x^{-1}\eta_{xx})_x| \{|(\bar{\rho}^{\gamma}\eta_x^{-\gamma})_x||\eta_{xx}|+|[\bar{\rho}(\eta_x^{-\gamma}-1)]_x|+\e^2|(\bar{\rho}v_t)_x|\}dxdt\\
&\le C\mathfrak{E}(0) + \delta \int_0^t\int_I\bar{\rho}^{4\gamma-3+3\varsigma} \eta_x^{-\gamma}|(\eta_x^{-1}\eta_{xx})_x|^2dxdt\\
&\ \   +C(\delta)\bar{Q}^{\gamma+2}\int_0^t\int_I (\bar{\rho}^{4\gamma-3+3\varsigma}\eta_{xx}^4
+\bar{\rho}^{2\gamma-1+3\varsigma}\eta_{xx}^2)dxdt\\
&\ \   +C(\delta)\bar{Q}^{\gamma+2}\int_0^t\int_I [\bar{\rho}^{1+3\varsigma} (\eta_x-1)^2+\e^4\bar{\rho}^{1+3\varsigma}v_t^2
+\e^4\bar{\rho}^{2\gamma-1+3\varsigma}v_{tx}^2\}dxdt,
\end{split}
\eee
with the help of \ef{11} and the fact that $\bar\rho\le 1$.
Note that
\bee\begin{split}
LHS\ge & \bar Q^{-2}\e^2\int_I\bar{\rho}^{3\gamma-3+3\varsigma}\eta_{xxx}^2dx+\gamma\bar Q^{-\gamma-2}\int_0^t\int_I\bar{\rho}^{4\gamma-3+3\varsigma}\eta_{xxx}^2dxdt\\
& -C\bar Q^{2}\e^2\int_I\bar{\rho}^{3\gamma-3+3\varsigma}\eta_{xx}^4dx-C\bar{Q}^{\gamma+4}\int_0^t\int_I\bar{\rho}^{4\gamma-3+3\varsigma}\eta_{xx}^4dxdt,
\end{split}\eee
with
\be\label{65-1}\begin{split}
 \e^2\int_I &\bar{\rho}^{3\gamma-3+3\varsigma}\eta_{xx}^4dx  \lesssim  \e^2\|\bar{\rho}^{2\gamma-2+2\varsigma}\eta_{xx}^2\|_{L^\infty} \int_I\bar{\rho}^{\gamma-1+\varsigma}\eta_{xx}^2dx \\
 \lesssim &  \exp(C\bar Q^2)  \mathfrak{E}(0)\e^2\Big|\int_x^1(\bar{\rho}^{2\gamma-2+2\varsigma}\eta_{yy}^2)_ydy\Big| \\
\lesssim &  \exp(C\bar Q^2)  \mathfrak{E}(0) \e^2 \int_x^1\lt(\bar{\rho}^{\gamma-1+\varsigma}\eta_{yy}^2+\bar{\rho}^{2\gamma-2+2\varsigma}|\eta_{yy}||\eta_{yyy}|\rt)dy \\
\lesssim & \delta\e^2\int_I\bar{\rho}^{3\gamma-3+3\varsigma}\eta_{xxx}^2dx+C(\delta) \exp(C\bar Q^2)  (\mathfrak{E}(0)+\mathfrak{E}^2(0))\e^2\int_I\bar{\rho}^{\gamma-1+\varsigma}\eta_{xx}^2dx,
\end{split}\ee
and similarly
\bee\begin{split}
  \int_0^t\int_I &\bar{\rho}^{4\gamma-3+3\varsigma}\eta_{xx}^4dxdt\\
  \lesssim &
  \delta\int_0^t\int_I\bar{\rho}^{4\gamma-3+3\varsigma}\eta_{xxx}^2dxdt
  +C(\delta)\exp(C\bar Q^2)  (\mathfrak{E}(0)+\mathfrak{E}^2(0))\int_0^t\int_I\bar{\rho}^{2\gamma-1+\varsigma}\eta_{xx}^2dxdt,
\end{split}\eee
where \ef{11} and \ef{41} have been used in \ef{65-1}.
Thus we derive
\ef{64} by using \ef{30B}, \ef{24}, \ef{41},  Lemma \ref{lem5} and the above inequalities. \hfill$\Box$.

\begin{lem}\label{lem11}(i) For $\gamma\in(1,3)$, we have
$$
\e^8\lt\|\lt(\frac{v_x}{\eta_x}\rt)_{xx}(\cdot,t)\rt\|_{L^2}^2\lesssim \exp(C\bar Q^2)  (\mathfrak{E}(0)+ \mathfrak{E}^7(0)),\quad t\ge 0.
$$
(ii) For $\gamma\in(1,\frac3 2)$, we obtain
\bee
\e^8\|v_{xxx}(\cdot,t)\|_{L^2}^2\lesssim \exp(C\bar Q^2)  (\mathfrak{E}(0)+ \mathfrak{E}^8(0)),\quad t\ge 0,
\eee
\bee
\e^8\|\eta_{xxx}(\cdot,t)\|_{L^2}^2\lesssim  (1+t)^{\frac{3\gamma-3+\varsigma}{\gamma}} \exp(C\bar Q^2)  (\mathfrak{E}(0)+ \mathfrak{E}^{10}(0)),\quad t\ge 0,
\eee
where $\varsigma$ is an arbitrary small positive constant.
\end{lem}
{\bf Proof.} (i) Differentiating \ef{8}  with respect to the spatial variable, we get
\bee\begin{split}
 \lt(\frac{v_x}{\eta_x}\rt)_{xx}= \bar{\rho}v_{tx}+& \bar{\rho}_xv_t
 +\frac{1}{\e^2}(\bar{\rho}^{\gamma})_{xx}(\eta_x^{-\gamma}-1)-\frac{2\gamma}{\e^2} (\bar{\rho}^{\gamma})_x\eta_x^{\gamma-1}\eta_{xx}\\
& -\frac{\gamma}{\e^2}\bar{\rho}^{\gamma}\lt[\eta_x^{-\gamma-1}\eta_{xxx}-(\gamma+1)\eta_x^{-\gamma-2}\eta_{xx}^2\rt].
\end{split}
\eee
Note that for $\gamma<3$, we have $\bar\rho^{2\gamma}\lesssim \bar\rho^{3\gamma-3+3\varsigma}$, where $\varsigma$ is a sufficiently small positive constant.
It follows that, for $\gamma\in (1,3)$ and suitably small $\varsigma>0$,
\bee\begin{split}
&\e^8\lt\|\lt(\frac{v_x}{\eta_x}\rt)_{xx}\rt\|_{L^2}^2\\
\lesssim & \int_I(\e^8\bar{\rho}^2v_{tx}^2+\e^8\bar{\rho}^{4-2\gamma}v_t^2+\bar{\rho}^{4-2\gamma}(\eta_x-1)^2+\bar{\rho}^2\eta_{xx}^2)dx +\e^2\int_I\bar{\rho}^{2\gamma}(\eta_{xxx}^2+\eta_{xx}^4)dx\\
\lesssim & \int_I[\e^8\bar{\rho}^2(v_{tx}^2+v_t^2)+\bar{\rho}^{1-\gamma+\varsigma}(\eta_x-1)^2
+\bar{\rho}^{\gamma-1+\varsigma}\eta_{xx}^2]dx+\e^2\int_I\bar{\rho}^{3\gamma-3+ \varsigma} \eta_{xxx}^2 dx\\
&+\e^2\int_I\bar{\rho}^{3\gamma-3+3\varsigma}   \eta_{xx}^4 dx\\
\lesssim & \exp(C\bar Q^2)  (\mathfrak{E}(0)+ \mathfrak{E}^7(0)),
\end{split}\eee
by using \ef{30B}, \ef{24}, \ef{41}, \ef{57}, \ef{64}, \ef{65-1} and Lemma \ref{lem33}.

(ii) Applying $\partial_x$ to \ef{50}, we find
\be\label{66}\begin{split}
  v_{xxx}=& \lt[\frac{\bar{\rho}^{\gamma}}{\e^2}(\eta_x^{-\gamma}-1)
  -\int_x^1\bar{\rho}v_tdy-\frac{\gamma}{\e^2}\bar{\rho}^{\gamma}\eta_x^{-\gamma}
  \rt]\eta_{xxx}
   \\
& +\lt[- \frac{2\bar{\rho}}{\e^2} (\eta_x^{-\gamma}-1)+\frac{\gamma(\gamma-1)}{\e^2}\bar{\rho}^{\gamma}\eta_x^{-\gamma-1}\eta_{xx}
+2\bar{\rho}v_t+\frac{2\gamma}{\e^2} \bar{\rho} \eta_x^{-\gamma}\rt]\eta_{xx}\\
& + \lt[-\frac{\bar{\rho}_x}{\e^2}(\eta_x^{-\gamma}-1) + \bar{\rho}_xv_t+\bar{\rho}v_{tx}\rt]\eta_x.
\end{split}\ee
Thus,
for $\gamma\in(1,\frac3 2)$ and $\varsigma>0$ sufficiently small, we have
\bee\begin{split}
 &\e^8\|v_{xxx}(\cdot,t)\|_{L^2}^2\\
\lesssim & \exp(C\bar Q^2)   \mathfrak{E}(0) \Big\{\e^2\int_I\bar{\rho}^{\gamma}\eta_{xxx}^2dx+\int_I\bar{\rho}^{2\cdot\min(1,2-\gamma)}\eta_{xx}^2dx \\
& + \e^2\int_I\bar{\rho}^{2\gamma}\eta_{xx}^4dx+\int_I\bar{\rho}^{4-2\gamma}[(\eta_x-1)^2+
\e^4v_t^2]dx+\e^8\int_I\bar{\rho}^2v_{tx}^2dx\Big\}\\
\lesssim & \exp(C\bar Q^2)   \mathfrak{E}(0) \lt\{\int_I\lt( \e^2\bar{\rho}^{3\gamma-3+\varsigma}\eta_{xxx}^2+ \bar{\rho}^{\gamma-1+\varsigma}\eta_{xx}^2+ (\eta_x-1)^2+ \e^8\bar\rho^2v_{tx}^2\rt)dx \rt. \\
& \lt.  +\e^2\int_I\bar{\rho}^{3\gamma-3+3\varsigma}  \eta_{xx}^4 dx\rt\}\\
\lesssim  & \exp(C\bar Q^2)  (\mathfrak{E}(0)+ \mathfrak{E}^8(0)),
\end{split}
\eee
by applying \ef{30B}, \ef{24}, \ef{41}, \ef{57}, \ef{64}, \ef{65-1} and Lemma \ref{lem33}. Note that we have also used the fact that for $\gamma<\frac 3 2$ and suitably small $\varsigma>0$, $\bar\rho^\gamma\lesssim \bar\rho^{3\gamma-3+\varsigma}$, $\bar{\rho}^{2\cdot\min(1,2-\gamma)}\lesssim \bar{\rho}^{\gamma-1+\varsigma}$, $\bar\rho^{2\gamma}\lesssim \bar\rho^{3\gamma-3+3\varsigma}$  and $\bar\rho^{4-2\gamma}\lesssim \bar\rho^{1+\varsigma}\lesssim 1$ here and thereafter.

Next, multiplying \ef{65} by $\e^8(\eta_x^{-1}\eta_{xx})_x$ and integrating over $I\times(0,t)$, we discover
\bee\begin{split}
 \e^8 \int_I\eta_{xxx}^2 & dx(t)+\e^6\int_0^t\int_I\bar{\rho}^{\gamma}\eta_{xxx}^2dxdt\\
\lesssim & \mathfrak{E}(0)+\e^8\int_I\eta_{xx}^4dx(t)+\e^6\int_0^t\int_I\bar{\rho}^{\gamma}\eta_{xx}^4dxdt
+\int_0^t\int_I\bar{\rho}^{2-\gamma}(\e^6\eta_{xx}^2+\e^8v_{tx}^2)dxdt \\
& + \int_0^t\int_I\bar{\rho}^{4-2\gamma}[\e^6(\eta_x-1)^2+\e^8v_t^2]dxdt  \\
\lesssim & \mathfrak{E}(0)+\e^8\int_I\eta_{xx}^2dx(t)\|\eta_{xx}^2\|_{L^\infty}(t)
+\e^2\sup_{t}\int_I\eta_{xx}^2dx(t)\e^4\int_0^t\|\bar{\rho}^{\gamma}\eta_{xx}^2\|_{L^\infty}dt\\
& +\e^8\int_0^t\lt(\int_I\bar{\rho}^{\gamma-1+\varsigma}\eta_{xx}^2dx\rt)^{\frac{3\gamma-3+\varsigma}{\gamma}}\lt(\int_I\bar{\rho}^{2\gamma-1+\varsigma}\eta_{xx}^2dx\rt)^{\frac{3-2\gamma-\varsigma}{\gamma}}dt \\
& +\int_0^t\int_I(\e^8v_{tx}^2+\e^8v_t^2+\bar{\rho}^{1+\varsigma}(\eta_x-1)^2)dxdt\\
\le & (1+t)^{\frac{3\gamma-3}{\gamma}+2\varsigma}\exp\{C\bar Q^2\}\lt(\mathfrak{E}(0)+\mathfrak{E}^{10}(0)\rt),
\end{split}\eee
by observing that
\bee\begin{split}
\e^8\int_I\eta_{xx}^2  dx\|\eta_{xx}^2\|_{L^\infty} \lesssim & \e^8\int_I\eta_{xx}^2dx\Big|\int_0^x\eta_{yy}\eta_{yyy}dy\Big| \\
\lesssim & \delta\e^8\int_I\eta_{xxx}^2dx+C(\delta)\lt(\e^2\int_I\eta_{xx}^2dx\rt)^3 \\
\lesssim & \delta\e^8\int_I\eta_{xxx}^2dx+C(\delta) (1+t)^{\frac{3(\gamma-1)}{\gamma}+3\varsigma}\exp(C\bar Q^2)  (\mathfrak{E}(0)+ \mathfrak{E}^9(0)),
\end{split}\eee
and
\bee\begin{split}
\e^2\sup_{t}\int_I&\eta_{xx}^2dx(t)\e^4\int_0^t\|\bar{\rho}^{\gamma}\eta_{xx}^2\|_{L^\infty}dt \\ \lesssim & \e^2\sup_{t}\int_I\eta_{xx}^2dx(t)\e^4\int_0^t\Big|\int_0^x(\bar{\rho}^{\gamma}\eta_{yy}^2)_ydy\Big|dt\\
\lesssim & \e^2\sup_{t}\int_I\eta_{xx}^2dx(t)\e^4\int_0^t\int_I(\bar{\rho}\eta_{xx}^2+\bar{\rho}^{2\gamma-1}\eta_{xxx}^2)dxdt\\
\lesssim & \e^2\sup_{t}\int_I\eta_{xx}^2dx(t)\lt(\e^4 \int_0^t\lt(\int_I\bar{\rho}^{\gamma-1+\varsigma}\eta_{xx}^2dx\rt)^{\frac{2\gamma-2+\varsigma}{\gamma}}
\lt(\int_I\bar{\rho}^{2\gamma-1+\varsigma}\eta_{xx}^2dx\rt)^{\frac{2-\gamma-\varsigma}{\gamma}}dt\rt.\\
& \lt. + \e^4 \int_0^t\lt(\int_I\bar{\rho}^{4\gamma-3+\varsigma}\eta_{xxx}^2dx\rt)^{\frac{\gamma-1}{3\gamma-3+\varsigma}}
\lt(\int_I\bar{\rho}^{\gamma}\eta_{xxx}^2dx\rt)^{\frac{2\gamma-2-\varsigma}{3\gamma-3+\varsigma}}dt\rt)\\
\lesssim & \delta \e^6\int_0^t\int_I\bar{\rho}^{\gamma}\eta_{xxx}^2dxdt+ C(\delta)(1+t)^{\frac{3\gamma-3+\varsigma}{\gamma}}(\mathfrak{E}(0)+ \mathfrak{E}^{10}(0)).
\end{split}\eee
Here \ef{30B}, \ef{24}, \ef{41}, \ef{49-1}, \ef{64} and Lemma \ref{hardy2} are utilized.
Therefore this lemma is shown.\hfill$\Box$
\vskip 3mm

\noindent\textbf{Proof of Theorem 2.2 (ii).} Let $\gamma\in (1,3/2)$. The energy estimates in Lemmas \ref{lem9}-\ref{lem11}, and the regularity proved in part (i) of this theorem, can give the high-order regularity in \ef{hregu}, by applying the refined Hardy inequality in Lemma \ref{lem33} with $w=v_{tt}$. The energy estimates in \ef{e2} are obtained by combining \ef{e}, \ef{57} and \ef{64}, and the decay estimates \ef{decay2} follow from \ef{57}.

Note that for each $\epsilon$, the energy estimates in \ef{e2} are global in the time variable. Moreover, since $\gamma\in (1,3/2)$, we find that the spatial derivative of each term in $\ef{6}_1$, namely $(\bar\rho v_t)_x$, $\e^{-2}(\bar\rho^\gamma\eta_x^{-\gamma})_{xx}=\e^{-2}[(\bar\rho^{\gamma})_{xx} \eta_x^{-\gamma} +2(\bar\rho^\gamma)_x (\eta_x^{-\gamma})_x + \bar\rho^\gamma (\eta_x^{-\gamma})_{xx}]$, $-\e^{-2}(\bar\rho g)_x$, $(v_x\eta_x^{-1})_{xx}=(2\eta_x^{-3}\eta_{xx}^2-\eta_x^{-2}\eta_{xxx})v_{x}
-2\eta_x^{-2}\eta_{xx}v_{xx} + \eta_x^{-1} v_{xxx}$, is bounded in $L^\infty(0,T;L^2(I))$, and the boundary terms in $\ef{6}_2$, $v$ and $v_x$ are bounded in $L^\infty(0,T;H^1(\partial I))$ by the trace theorem, for fixed $\epsilon$. Together with the estimates shown in the part (i) and the local existence result in Theorem 2.1, we can derive the global existence and uniqueness of the classical solution to \ef{6}.

\hfill$\Box$

\section{Singular limit}

{\bf Proof of Theorem \ref{thm3}.} Let $(\eta^\e,v^\e)$ be the solution obtained in Theorem \ref{thm2}. Due to \ef{e}, we have
\be\label{C19}
\int_0^t\int_0^{1-1/n} (\eta_x^\epsilon -1)^2 dx ds \le C(n)\epsilon^2,\quad \forall \; n\in\mathbb{N},
\ee
which gives
\bee
\eta_x^\epsilon(x,t)\to 1\quad{\rm as}\quad \epsilon\to 0,\; a.e.\; (x,s)\in (0,1-1/n)\times(0,t).
\eee
Since $(0,1)=\cup_{n=1}^\infty (0, 1-1/n)$, we find
\be\label{C20}
\eta_x^\epsilon(x,t)\to 1\quad{\rm as}\quad \epsilon\to 0,\; a.e.\; (x,s)\in I\times(0,t), t\in (0,T].
\ee
On the other hand, we derive from \ef{e} that, for any fixed $p\ge 1$, $\{(\eta_x^\epsilon)^p\}_{\epsilon\in (0,1]}$ and $\{(\eta_x^\epsilon)^{-p}\}_{\epsilon\in (0,1]}$ are bounded in $L^\infty(0,T; L^\infty(I))$; moreover, from \ef{e1},
$\{\partial_t[(\eta_x^\epsilon)^p]\}_{\epsilon\in (0,1]}$ and $\{\partial_t[(\eta_x^\epsilon)^{-p}]\}_{\epsilon\in (0,1]}$ are bounded in $L^2(0,T; L^2(I))$, since $\partial_t[(\eta_x^\epsilon)^p]=p(\eta_x^\epsilon)^{p-1} v_x^\epsilon$ and $\partial_t[(\eta_x^\epsilon)^{-p}]=-p(\eta_x^\epsilon)^{-p-1} v_x^\epsilon$. Thus, together with \ef{C20}, we find
\be\label{C21}
(\eta_x^\epsilon)^p\to 1,\; (\eta_x^\epsilon)^{-p}\to 1\q {\rm strongly} \; {\rm in} \; C([0,T],L^2(I)),\q\epsilon\to 0,
\ee
by using Lemma \ref{compact}.
It follows that
\be\label{C22}
f^\epsilon=\rho^\epsilon(\eta^\epsilon)=\bar\rho/\eta_x^\epsilon \to \bar\rho \q   {\rm strongly} \; {\rm in} \; C([0,T],L^2(I)),\q\epsilon\to 0.
\ee
Note that from \ef{eta0def}, we find $\eta_0^\epsilon(0)=0$, which gives that
\bee
\eta^\epsilon(0,t)=\eta_0^\epsilon(0)+\int_0^t v^\epsilon(0,s)ds=0.
\eee
Thus from \ef{C21} and the above equality, we obtain
\be\label{C21p}
\eta^\epsilon(x,t)=  \int_0^x\eta_y(y,t)dy\to x \quad {\rm strongly} \; {\rm in}\quad C([0,T],L^2(I)),
\ee
and the free boundary
\be\label{C22p}
\eta^\epsilon(1,t)=  \int_0^1\eta_x(1,t)dx\to 1,\quad a.e.\quad t\in[0,T].
\ee
Next, from \ef{e}, \ef{e1} and Lemma \ref{lem33}, there exists a function $v(x,t)$, such that
\be\label{C23}
v^\epsilon\rightharpoonup v,\; v_x^\epsilon\rightharpoonup v_x\; {\rm weakly} \; {\rm in} \; L^2( 0,T; L^2(I)), \;\epsilon\to 0.
\ee
Now it suffices to show $v=0$ and $v_x=0$ a.e. $(x,t)\in I\times(0,t)$. Indeed, from \ef{419}, the following identity holds:
\be\label{C24}
-\int_0^t\int_I f^\epsilon \varphi_t  dxdt+ \int_0^t\int_I f^\epsilon (\eta_x^\epsilon)^{-1} v_x^\epsilon\varphi  dxdt=\int_I\bar\rho (\eta_{0x}^\epsilon)^{-1} \varphi(x,0) dx,
\ee
for any $\varphi\in C_c^1([0,t),C_c^\infty(I))$. Letting $\epsilon \to 0,$ we find
\bee
\begin{split}
-\int_0^t &\int_I f^\epsilon \varphi_t  dxdt=- \int_0^t\int_I \bar\rho (\eta_x^\epsilon)^{-1} \varphi_t  dxdt\\
 & \to -\int_0^t\int_I \bar\rho \varphi_t dxdt=\int_I \bar\rho(x) \varphi(x,0)dx,\quad \epsilon\to 0,
\end{split}
\eee
and
\bee
\begin{split}
   \int_I \bar\rho (\eta_{0x}^\epsilon)^{-1} \varphi(x,0)  dx  \to \int_I \bar\rho(x) \varphi(x,0)dx,\quad \epsilon\to 0,
\end{split}
\eee
thus
\be\label{C25}
\int_0^t\int_I f^\epsilon (\eta_x^\epsilon)^{-1} v_x^\epsilon\varphi  dxdt \to 0,\quad \epsilon\to 0.
\ee
On the other hand, we apply \ef{C21} and \ef{C23} to get
\be\label{C26}
\int_0^t\int_I f^\epsilon (\eta_x^\epsilon)^{-1} v_x^\epsilon\varphi  dxdt=\int_0^t\int_I\bar\rho \varphi (\eta_x^\epsilon)^{-2} v_x^\epsilon   dxdt \to \int_0^t\int_I\bar\rho v_x \varphi dxdt,\quad \epsilon\to 0.
\ee
Therefore, from \ef{C25} and \ef{C26}, we infer
\bee
\bar\rho v_x=0, \quad a.e.\; (x,t)\in I\times(0,t),\;t\in (0,T],
\eee
which yields
\be\label{C27}
 v_x=0, \quad a.e.\; (x,t)\in I\times(0,t),\;t\in (0,T],
\ee
since $\bar\rho (x)$ only vanishes at $x=1$. Invoking of \ef{C27} and the boundary condition $v(0,t)=0$, one shows
\be\label{C28}
 v =0, \quad a.e.\; (x,t)\in I\times(0,t),\;t\in (0,T].
\ee
Therefore, \ef{conv1} follows from \ef{C21}, \ef{C21p}, \ef{C23}, \ef{C27} and \ef{C28}; \ef{conv2} comes from \ef{C22}, \ef{conv1} and the fact that $u^\epsilon(\eta^\epsilon(x,t))=v^\epsilon(x,t)$; and \ef{conv3} is shown by \ef{C22p}.
\hfill$\Box$

\section{Local existence}
In this section, we show the local existence and uniqueness of the classical solution to \ef{6} by the finite difference scheme, which is also applied in \cite{LiXY,LXZ}.
\\
{\bf Proof of Theorem 2.1.} Since we aim to prove the local in time existence of solutions, without the loss of generality, we take $\epsilon=1$. Let $N$ be a positive integer and $h=1/N$. Then for $i=0,1,\cdots,N$, we set the end points $x_i=ih$, and
$$\bar\rho^i=\bar\rho(x_i),\q \eta^i(t)=\eta(x_i,t), \q v^i(t)=v(x_i,t).$$
In particular, $\bar\rho^i >0$ for $i=0,1,\cdots,N-1$ and $\bar\rho^N=0.$
  Since the aim is to prove the existence of classical solutions, we replace  $\ef{6}_1$ by an equivalent equation \ef{8}, and approximate the resulting system \ef{3}-\ef{8}-$\ef{6}_{2,3}$  by the following Cauchy problem of ordinary differential  equations for $(\eta^i,v^i)$, $i=1,\cdots,N-1$, for $t>0$:
\be\label{A1}\left\{
\begin{split}
&(\eta^i) '=v^i,\\
&\bar\rho^i (v^i) ' +\frac{1}{h}\lt\{(\bar\rho^{i+1})^\gamma \lt[\frac{h^\gamma}{(\eta^{i+1}-\eta^i)^\gamma}-1\rt]-(\bar\rho^{i})^\gamma \lt[\frac{h^\gamma}{(\eta^{i}-\eta^{i-1})^\gamma}-1\rt]\rt\}\\
&\qquad =\frac{1}{h}\lt(\frac{v^{i+1}-v^i}{ \eta^{i+1}-\eta^i}
-\frac{v^{i}-v^{i-1}}{ \eta^{i}-\eta^{i-1} }\rt),\\
&\eta^i(0)=\eta_0(x_i),\; v^i(0)=v_0(x_i).
\end{split}\right.
\ee
Moreover, from the boundary conditions in $\ef{6}_2$, we set
  \be\label{A2}
  v^0(t)=0,\q v^N(t)=v^{N-1}(t),\q t\ge0,
  \ee
  with
  \be\label{A2p}
\partial_t v^0(0):=\partial_t v^0(0_+)=0,\; \partial_t (v^N-v^{N-1})(0):=\partial_t (v^N-v^{N-1})(0_+)=0;
\ee
in addition,
 since $\eta_0(0)=0$ and $\eta_t=v$, we set
 \be\label{A3}
 \eta^0(t)=0,\q \eta^N(t)=\eta^{N-1}(t)+\eta_0(x_N)-\eta_0(x_{N-1}),\q t\ge0.
 \ee
For the sake of completeness in the following definition of $\mathfrak{E}_N(t)$, we define additionally
\be\label{A3p}
v^{-1}(t)=0,\;\eta^{-1}(t)=0,\; t\ge 0,
 \ee
 \be\label{A3pp}
v^{N+1}(t)=v^{N}(t),\;\eta^{N+1}(t)=\eta^{N}(t)+\eta_0(x_{N+1})-\eta_0(x_{N}),\; t\ge0,
 \ee
 with $\partial_t v^{-1}(0)=0$ and  $\partial_t (v^{N+1}-v^N)(0)=0$.

Define the functional $\mathfrak{E}_N(t)$ for the approximate problem \ef{A1} by
\be\label{A4}
\begin{split}
\mathfrak{E}_N(t):=& \max_{1\le i\le N} \lt\{\lt|\frac{\eta^i(t)-\eta^{i-1}(t)}{h}\rt|^2 + \lt|\frac{h}{\eta^i(t)-\eta^{i-1}(t)}\rt|^2 +\lt|\frac{v^i(t)-v^{i-1}(t)}{h}\rt|^2 \rt\}\\
&+h\sum_{i=1}^{N-1}\bar\rho^i(|v^i |^2+|(v^i)'|^2 + |(v^i)''|^2) +  h\sum_{i=1}^{N-1}  \lt|\frac{v^{i+1}(t)-2v^i(t)+v^{i-1}(t)}{h^2}\rt|^2\\
& +  h\sum_{i=1}^{N-1}  \lt|\frac{\eta^{i+1}(t)-2\eta^i(t)+\eta^{i-1}(t)}{h^2}\rt|^2\\
 &+ h\sum_{i=1}^{N}  \lt|\frac{v^{i+1}(t)-3v^i(t)+3v^{i-1}(t)-v^{i-2}(t)}{h^3}\rt|^2\\
&+ h\sum_{i=1}^{N}  \lt|\frac{\eta^{i+1}(t)-3\eta^i(t)+3\eta^{i-1}(t)-\eta^{i-2}(t)}{h^3}\rt|^2.
\end{split}
\ee
Note that the definition of $ {\mathfrak{E}}_N(t)$ is not exactly the discrete version of $\mathfrak{E}(t)$ in \ef{Et}, since the uniform in time estimates are not required. For instance, for $\alpha>0$, we can estimate the weighted terms
$$\int_I \bar\rho^{1-\gamma+\alpha}(\eta_x-1)^2 dx (t) + \int_0^t\int_I \bar\rho^{1+\alpha} (\eta_x-1)^2 dx dt \le C(t) (\|\eta_x \|_{L^\infty_{x,t}} +1)\le C(t),$$
provided that $\|\eta_x \|_{L^\infty_{x,t}}$ is bounded.
Instead, the continuous version of $\mathfrak{E}_N(t)$ is defined by
\be\label{BarEt}
\begin{split}
\bar{\mathfrak{E}}(t):=& \|\eta_x(\cdot,t)\|_{L^\infty}^2 + \|\eta_x^{-1}(\cdot,t)\|_{L^\infty}^2 + \|v_x(\cdot,t)\|_{L^\infty}^2+ \int_I \bar\rho(v^2 + v_t^2 +v_{tt}^2) dx (t)\\
 &+ \int_I (v_{xx}^2 + \eta_{xx}^2+v_{xxx}^2 + \eta_{xxx}^2)dx (t).
\end{split}
\ee
 For sufficiently large $N$, we have $\mathfrak{E}_N(0)\le 2\bar{\mathfrak{E}}(0)$. Note that $(v^i)'(0)$, $i=1,...,N-1$, can be expressed by taking $t=0$ in  $\ef{A1}_2$, that is,
 \bee
\begin{split}
 (v^i) '(0)= & \frac{1 }{\bar\rho^i h}\lt(\frac{v^{i+1}(0)-v^i(0)}{ \eta^{i+1}(0)-\eta^i(0)}
-\frac{v^{i}(0)-v^{i-1}(0)}{ \eta^{i}(0)-\eta^{i-1}(0) }\rt.\\
&\lt.-(\bar\rho^{i+1})^\gamma \lt[\frac{h^\gamma}{(\eta^{i+1}(0)-\eta^i(0))^\gamma}-1\rt]+ (\bar\rho^{i})^\gamma \lt[\frac{h^\gamma}{(\eta^{i}(0)-\eta^{i-1}(0))^\gamma}-1\rt]\rt),
\end{split}
\eee
while   $(v^{N}) '(0)=(v^{N-1}) '(0)$ due to the assumption \ef{A2p}. Thus by $\bar\rho^N=0$ and the properties of integrals, we have
 $$h\sum_{i=1}^{N-1}\bar\rho^i |(v^i)' (0)|^2 =h\sum_{i=1}^{N}\bar\rho^i |(v^i)' (0)|^2 \le 2 \int_I \bar\rho v_t^2(x,0)dx,$$
 for sufficiently large $N$. The remained terms in $\mathfrak{E}_N(0)$ can be bounded by the similar way.
\begin{lem}\label{uniform} {\rm (Uniform estimates in $N$)} There exists a positive constant $T^*$ independent of $N$, such that for any $t\in [0,T^*]$,
\be\label{AA1}
\mathfrak{E}_N(t)+ \frac{1}{h}\int_0^t\sum_{i=1}^{N-1}  (|(v^i)'-(v^{i-1})'|^2+ |(v^i)''-(v^{i-1})''|^2)(s)ds\le 2\bar C \mathfrak{E}_N(0) ,
\ee
where the positive constant $\bar{C}$ is independent of $N$ and $t\in [0,T^*]$.
\end{lem}
{\bf Proof.} By the general theory for ordinary differential equations,  \ef{A1} admits a  local  solution  where the maximal time of existence is denoted as $T_N$. Observe that $t\sup_{s\in[0,t]}\mathfrak{E}_N(s)$ is a continuous function of $t$ and choose $T_N'\le T_N$ to satisfy
\be\label{A6}
T_N'\sup_{s\in[0,T_N']}\mathfrak{E}_N(s)\le ( 2 D_0)^{-1}\le   D_0 /2.
\ee
Then from $\ef{A1}_1$, for $i=1,\cdots,N-1$ and $t\in [0,T_N'),$
\be\label{A7}
\begin{split}
 \frac{\eta^i(t)-\eta^{i-1}(t)}{h} &\le \lt|\frac{\eta^i(0)-\eta^{i-1}(0)}{h}\rt| + \int_0^t \lt|\frac{v^i(t)-v^{i-1}(t)}{h}\rt|\\
&\le \|\eta_{0x}\|_{L^\infty}+ t\sup_{s\in[0,t]}\mathfrak{E}_N(s)\le \frac 3 2 D_0.
\end{split}
\ee
Similarly,
 \be\label{A8}
\begin{split}
 \frac{\eta^i(t)-\eta^{i-1}(t)}{h} &\ge  \inf_{x\in I} \eta_{0x}(x)- T_N'\sup_{s\in[0,T_N']}\mathfrak{E}_N(s) \ge   \frac 1 2 D_0^{-1}.
\end{split}
\ee
In what follows, we only need to prove that $T_N'$ has a positive lower bound which is independent of $N$.
To this end , we should derive an inequality in the form:
\be\label{A9}
\mathfrak{E}_N(t)\le \bar C \mathfrak{E}_N(0) + \bar C \int_0^t (\mathfrak{E}_N(s) + \mathfrak{E}_N^2(s))ds,
\ee
where the constant $\bar C>1$ is to be determined. Indeed, we can set $$T_N^*:=\sup\{t\in[0,T_N']| \mathfrak{E}_N(t)\le 2\bar C \mathfrak{E}_N(0) \}.$$
Clearly we have  $T_N^*>0$ since $\bar C>1.$ Suppose that $T_N^*<({2D_0 \bar C \mathfrak{E}_N(0)})^{-1}$ (Otherwise, $T_N^*=+\infty$ or $T_N^*$ is a positive constant satisfying  $T_N^*\ge ({2D_0 \bar C \mathfrak{E}_N(0)})^{-1}\ge ({4D_0 \bar C \mathfrak{E} (0)})^{-1}$, thus we can obtain a positive lower bound of existence independent of $N$). Due to the definition of $T_N^*$, we find that $T_N^*$ also satisfies \ef{A6}, moreover,
\bee
\begin{split}
  2\bar C \mathfrak{E}_N(0)=\mathfrak{E}_N(T_N^*)&\le \bar C \mathfrak{E}_N(0) + \bar C \int_0^{T_N^*} (\mathfrak{E}_N(s) + \mathfrak{E}_N^2(s))ds\\
  &\le  \bar{C} \mathfrak{E}_N(0) + 2 \bar C T_N^* \mathfrak{E}_N(0)(1+ 2\bar C\mathfrak{E}_N(0)).
\end{split}
\eee
It follows that
\bee
\begin{split}
    T_N^* \ge [2 \bar C  (1+ 2\bar C\mathfrak{E}_N(0))]^{-1}\ge [2 \bar C   (1+ 4\bar C\mathfrak{E} (0))]^{-1}=: T^*,
\end{split}
\eee
which gives \ef{AA1} immediately.

We are ready to show \ef{A9} term by term. Clearly, from $\ef{A1}_1$ and \ef{A3}, we have
\be\label{A10}
\max_{1\le i\le N}  \lt|\frac{\eta^i(t)-\eta^{i-1}(t)}{h}\rt|^2 \le \mathfrak{E}_N(0) +  \int_0^t \mathfrak{E}_N(s)ds,
\ee
where the calculations are similar to \ef{A7}.
Next, multiplying $\ef{A1}_2$ by $h v^i(t)$ and summing up the resulting equalities for $i=1,\cdots,N-1$, we get
\be\label{A11}
\begin{split}
 \frac{d}{dt}&\sum_{i=1}^{N-1}\lt\{\bar\rho^i \frac{|v^i|^2}{2} +  (\bar\rho^{i})^\gamma \lt[\frac{1}{\gamma-1}\lt(\frac{ \eta^{i}-\eta^{i-1}}{h}\rt)^{1-\gamma}+ \frac{ \eta^{i}-\eta^{i-1}}{h}-\frac{\gamma}{\gamma-1}\rt]\rt\} h\\
  &+ \sum_{i=1}^{N-1}\frac{(v^i-v^{i-1})^2}{\eta^i-\eta^{i-1}}=0,
\end{split}
\ee
which gives
\be\label{A12}
\begin{split}
   \sum_{i=1}^{N-1}&\lt[\bar\rho^i |v^i(t)|^2  +(\bar\rho^{i})^\gamma\lt(\frac{ \eta^{i}(t)-\eta^{i-1}(t)}{h}-1\rt)^2 \rt]h \\
     &+ \int_0^t\sum_{i=1}^{N-1} \frac{|v^i(s)-v^{i-1}(s)|^2}{h}ds\le C\mathfrak{E}_N(0),
\end{split}
\ee
by using
\ef{A4}, \ef{A7}, \ef{A8} and the observation that
\bee
\frac{1}{\gamma-1}\lt(\frac{ \eta^i-\eta^{i-1}}{h}\rt)^{1-\gamma}+ \frac{\eta^i-\eta^{i-1}}{h}-\frac{\gamma}{\gamma-1}\sim \frac \gamma 2 \lt(\frac{ \eta^i-\eta^{i-1}}{h}-1\rt)^2.
\eee
Next, we sum up $2(\eta^{n}-\eta^{n-1})^{-2}h^3\cdot\ef{A1}_2$ for $i=n$ to $N-1$ to get
\be\label{A13}
\frac{d}{dt}\lt(\frac{h^2}{(\eta^n-\eta^{n-1})^2}\rt)=  \frac{2h^2}{(\eta^n-\eta^{n-1})^2}   \lt(\sum_{i=n}^{N-1} \bar\rho^i (v^i) ' h-(\bar\rho^{n})^\gamma \lt[\frac{h^\gamma}{(\eta^{n}-\eta^{n-1})^\gamma}-1\rt]\rt).
\ee
Integrating \ef{A13}, then using  \ef{A3}, \ef{A7}, \ef{A8} and \ef{A12}, and Taylor's expansion $y^\gamma=1 + \gamma \xi^{\gamma-1}(y-1)$, with $\xi$ being between $y$ and 1,  we find
\be\label{A14}
\max_{1\le i\le N}\lt|\frac{h}{ \eta^n(t)-\eta^{n-1}(t) }\rt|^2 \le C \lt( \mathfrak{E}_N(0)+\int_0^t \mathfrak{E}_N^\frac 3 2(s) ds \rt).
\ee
Differentiating $\ef{A1}_2$ yields
\be\label{A15}
\begin{split} \bar\rho^i (v^i)'' &+\frac{1}{h}\lt\{(\bar\rho^{i+1})^\gamma \lt[\frac{h^\gamma}{(\eta^{i+1}-\eta^i)^\gamma}\rt] '-(\bar\rho^{i})^\gamma \lt[\frac{h^\gamma}{(\eta^{i}-\eta^{i-1})^\gamma}\rt] '\rt\}\\
&  =\frac{1}{h}\lt(\frac{v^{i+1}-v^i}{ \eta^{i+1}-\eta^i}
-\frac{v^{i}-v^{i-1}}{ \eta^{i}-\eta^{i-1} }\rt)' .
\end{split}
\ee
Next,  multiplying \ef{A15} by $h (v^i)'(t)$, summing up the resulting equalities for $i=1,\cdots,N-1$, and integrating in $[0,t]$, we derive
\be\label{A17}
\begin{split}
   \sum_{i=1}^{N-1}& \bar\rho^i \frac{|(v^i)'(t)|^2}{2} h  + \int_0^t\sum_{i=1}^{N-1} \frac{((v^i)'(s)-(v^{i-1}(s))')^2}{\eta^{i}(s)-\eta^{i-1}(s)}ds\\
   &= \sum_{i=1}^{N-1}  \bar\rho^i \frac{|(v^i)'(0)|^2}{2} h+ \int_0^t\sum_{i=1}^{N-1} \frac{ [ (v^i(s))' -(v^{i-1}(s) )' ]}{(\eta^{i}(s)-\eta^{i-1}(s))^2}(v^{i}(s)-v^{i-1}(s))^2 ds\\
   &\quad -\gamma\int_0^t\sum_{i=1}^{N-1}  \bar\rho^{i}\frac{h^\gamma(v^i(s)-v^{i-1}(s))}{(\eta^i(s)-\eta^{i-1}(s))^{\gamma+1}}
   [(v^i(s))'-(v^{i-1}(s))']ds\\
   &\le C\mathfrak{E}_N(0)+\frac{1}{2}\int_0^t\sum_{i=1}^{N-1} \frac{((v^i)'(s)-(v^{i-1}(s))')^2}{\eta^{i}(s)-\eta^{i-1}(s)}ds\\
    &\quad+ C\int_0^t\sum_{i=1}^{N-1} \lt(\frac{|v^{i}(s)-v^{i-1}(s)|^4}{|\eta^{i}(s)-\eta^{i-1}(s)|^3} +\frac{h^{2\gamma}|v^{i}(s)-v^{i-1}(s)|^2}{ |\eta^{i}(s)-\eta^{i-1}(s)|^{2\gamma+1}}\rt)  ds,
  \end{split}
\ee
where $\bar\rho^N=0$ and \ef{A2}-\ef{A3} are used.
Thus, in view of \ef{A7} and \ef{A8}, one shows
\be\label{A18}
\begin{split}
  \sum_{i=1}^{N-1}& \bar\rho^i |(v^i)'(t)|^2 h +  \int_0^t\sum_{i=1}^{N-1} \frac{((v^i)'(s)-(v^{i-1}(s))')^2}{h^2}hds\\
  &\le C\mathfrak{E}_N(0) + C\int_0^t (\mathfrak{E}_N(s) + \mathfrak{E}_N^2(s))ds.
\end{split}
\ee
Sum up $\ef{A1}_2$ for $i=j$ to $N-1$ to get
\be\label{A19}
\frac{v^j-v^{j-1}}{h}=\frac{\eta^j-\eta^{j-1}}{h}\lt((\bar\rho^{j})^\gamma \lt[\frac{h^\gamma}{(\eta^{j}-\eta^{j-1})^\gamma}-1\rt]-\sum_{i=j}^N \bar\rho^i (v^i) ' h\rt),
\ee
It yields
\be\label{A20}
\lt|\frac{v^j(t)-v^{j-1}(t)}{h}\rt|^2\le  C\mathfrak{E}_N(0) + C\int_0^t \mathfrak{E}_N^2(s)ds,\quad j=1,\cdot,N,
\ee
by applying \ef{A2}, \ef{A7}, \ef{A8}, \ef{A18} and Taylor's expansion
$$y^{-\gamma}=1 - \gamma \xi^{-\gamma-1}(y-1),$$
where $\xi$ is between $y$ and 1.
 Moreover, similarly as \ef{A18}, we can derive the estimates
\be\label{A21}
\begin{split}
  \sum_{i=1}^{N-1}& \bar\rho^i |(v^i)''(t)|^2 h +  \int_0^t\sum_{i=1}^{N-1} \frac{((v^i)''(s)-(v^{i-1}(s))'')^2}{h}ds\\
  &\le C\mathfrak{E}_N(0) + C\int_0^t \mathfrak{E}_N^2(s)ds,
\end{split}
\ee
 by differentiating \ef{A15} in $t$, integrating the product of the resulting equality and $(v^i)''(t)$, summing up in $i$, integrating in $t$ and using \ef{A2}-\ef{A3} and $\bar\rho^N=0$.

To get the estimates for higher order difference scheme, we rewrite $\ef{A1}_2$ as
\be\label{A22}
\begin{split}
\frac{1}{h}&\lt(\log\frac{\eta^{i+1}-\eta^i}{\eta^i-\eta^{i-1}}\rt)_t=  \bar\rho^i (v^i) '\\
& +\frac{1}{h}\lt\{(\bar\rho^{i+1})^\gamma \lt[\frac{h^\gamma}{(\eta^{i+1}-\eta^i)^\gamma}-1\rt]-(\bar\rho^{i})^\gamma \lt[\frac{h^\gamma}{(\eta^{i}-\eta^{i-1})^\gamma}-1\rt]\rt\}.
\end{split}
\ee
Next, we multiply \ef{A22} by  $\log\frac{\eta^{i+1}-\eta^i}{\eta^i-\eta^{i-1}}$, sum up with respect to $i$, then integrate in $t$ and using Taylor's expansions
$$\log y=\xi_1^{-1} (y-1),\quad y^{-\gamma}=z^{-\gamma}-\gamma \xi_2^{-\gamma-1}(y-z),$$
where $\xi_1$  is between $y$ and 1, and $\xi_2$ is between $y$ and $z$, to get
\be\label{A23}
\begin{split}
\sum_{i=1}^{N-1}&\frac{1}{h} \lt         ( \frac{\eta^{i+1}-\eta^i}{\eta^i-\eta^{i-1}}-1\rt)^2 \le C\mathfrak{E}_N(0) + C\int_0^t\sum_{i=1}^{N-1}   \bar\rho^i |(v^i) '|^2 h ds\\
 & + C\int_0^t\sum_{i=1}^{N-1}  \lt(  \lt|\frac {\eta^{i}-\eta^{i-1}}{h} -1\rt|^2 + \lt|\frac{\eta^{i+1}(t)-2\eta^i(t)+\eta^{i-1}(t)}{h^2}\rt|^2 \rt)h ds ,
\end{split}
\ee
which gives
\be\label{A24}
\begin{split}
\sum_{i=1}^{N-1} \lt|\frac{\eta^{i+1}(t)-2\eta^i(t)+\eta^{i-1}(t)}{h^2}\rt|^2 h
 \le C\mathfrak{E}_N(0) + C\int_0^t \mathfrak{E}_N(s) ds,
\end{split}
\ee
by using \ef{A3}, \ef{A3pp} and \ef{A7}.
In virtue of  \ef{A19}, we derive
\be\label{A25}
\begin{split}
\frac{v^{i+1} -2v^i +v^{i-1} }{h^2}=  &\frac{\eta^{i+1} -2\eta^i +\eta^{i-1} }{h^2}\lt[ (\bar\rho^{i+1})^\gamma\lt(\frac{h^\gamma}{(\eta^{i+1}-\eta^i)^\gamma}-1\rt)-\sum_{j=i+1}^N \bar\rho^j (v^j)'h\rt]\\
&+\frac{\eta^i-\eta^{i-1}}{h}
\lt[ \frac{(\bar\rho^{i+1})^\gamma-(\bar\rho^{i})^\gamma}{h}
 \lt(\frac{h^\gamma}{(\eta^{i+1}-\eta^i)^\gamma}-1\rt)\rt.\\
 &\left.  + \frac{(\bar\rho^{i})^\gamma}{h}\lt(\frac{h^\gamma}{(\eta^{i+1}-\eta^i)^\gamma}-
 \frac{h^\gamma}{(\eta^{i}-\eta^{i-1})^\gamma}\rt)+  \bar\rho^i (v^i)'\rt],
\end{split}
\ee
which yields
\be\label{A26}
\begin{split}
\sum_{i=1}^{N-1} \lt|\frac{v^{i+1} -2v^i +v^{i-1} }{h^2}\rt|^2 h
 \le C\mathfrak{E}_N(0) + C\int_0^t ( \mathfrak{E}_N(s) + \mathfrak{E}_N^2(s)) ds,
\end{split}
\ee
by using \ef{A2}, \ef{A3}, \ef{A3pp}, \ef{A7}, \ef{A8}, \ef{A24}, \ef{A10}, \ef{A14}, \ef{A18} and Taylor's expansion $y^{-\gamma}=z^{-\gamma}-\gamma\xi^{-\gamma-1} (y-z)$ for $\xi$ being between $y$ and $z$.
From \ef{A25} again, we have
\be\label{A27}
\begin{split}
 &\frac{v^{i+1}   -3v^i +3v^{i-1} -v^{i-2} }{h^3}\\
    = & \frac{\eta^{i+1} -3\eta^i +3\eta^{i-1} -\eta^{i-2} }{h^3} \lt[ (\bar\rho^{i+1})^\gamma\lt(\frac{h^\gamma}{(\eta^{i+1}-\eta^i)^\gamma}-1\rt)-\sum_{j=i+1}^N \bar\rho^j (v^j)'h\rt]\\
 &+ 2\frac{\eta^{i} -2\eta^{i-1} +\eta^{i-2} }{h^2} \frac 1 h  \lt[ (\bar\rho^{i+1})^\gamma\lt(\frac{h^\gamma}{(\eta^{i+1}-\eta^i)^\gamma}-1\rt)
 -(\bar\rho^{i})^\gamma\lt(\frac{h^\gamma}{(\eta^{i}-\eta^{i-1})^\gamma}-1\rt)+  \bar\rho^i (v^i)'h \rt] .\\
 &-\frac{\eta^{i-1}-\eta^{i-2}}{h}\frac{1}{h^2}\lt[ (\bar\rho^{i+1})^\gamma\lt(\frac{h^\gamma}{(\eta^{i+1}-\eta^i)^\gamma}-1\rt)-2 (\bar\rho^{i})^\gamma\lt(\frac{h^\gamma}{(\eta^{i }-\eta^{i-1})^\gamma}-1\rt) \rt.\\
&\quad \lt. +(\bar\rho^{i-1})^\gamma\lt(\frac{h^\gamma}{(\eta^{i-1}-\eta^{i-2})^\gamma}-1\rt)
+  (\bar\rho^{i} (v^{i})'-\bar\rho^{i-1} (v^{i-1}))'h\rt].
\end{split}
\ee
Observing that
\be\begin{split}
&\frac{\eta^{i+1} -3\eta^i +3\eta^{i-1} -\eta^{i-2} }{h^3}\\
&\qquad=\frac{\eta^{i+1}(0) -3\eta^i(0) +3\eta^{i-1}(0) -\eta^{i-2}(0) }{h^3} + \int_0^t \frac{v^{i+1}   -3v^i +3v^{i-1} -v^{i-2} }{h^3} ds,
\end{split}\ee
we derive
\be\label{A28}
\begin{split}
\sum_{i=1}^{N} &\lt(\lt|\frac{\eta^{i+1} -3\eta^i +3\eta^{i-1} -\eta^{i-2} }{h^3}\rt|^2+\lt|\frac{v^{i+1}   -3v^i +3v^{i-1} -v^{i-2} }{h^3}\rt|^2 \rt)h\\
 &\le C\mathfrak{E}_N(0) + C\int_0^t ( \mathfrak{E}_N(s) + \mathfrak{E}_N^2(s))  ds,
\end{split}
\ee
by applying \ef{A2}, \ef{A3}-\ef{A3pp}, \ef{A7}, \ef{A8}, \ef{A10}, \ef{A14}, \ef{A18}, \ef{A24}, \ef{A27} and the mean value theorem.
Note that, summing up \ef{A15} for $j=i,\cdots,N-1$ and multiplying the resulting equality by $h$, we find
\be \label{vtx0}
\begin{split}
\lt(\frac{v^i-v^{i-1}}{\eta^i-\eta^{i-1}}\rt)' &= -\sum_{j=i}^{N-1} \bar\rho^i (v^i)'' h + (\bar\rho^{i})^\gamma\lt[\frac{h^\gamma}{(\eta^i-\eta^{i-1})}\rt]'\\
&=-\sum_{j=i}^{N-1} \bar\rho^i (v^i)'' h -\gamma (\bar\rho^{i})^\gamma \frac{h^{\gamma+1}}{(\eta^i-\eta^{i-1})^{\gamma+1}} \frac{v^i-v^{i-1}}{h},
\end{split}
\ee
where the boundary conditions $\bar\rho^N=0$ and $v^N=v^{N-1}$ are used. It follows from \ef{A7}, \ef{A8}, \ef{A10}, \ef{A20}, \ef{A21} and \ef{vtx0} that
\be  \label{vtx}
\begin{split}
\lt| \frac{(v^i)'-(v^{i-1})'}{h} \rt|^2 &= \left| \frac{\eta^i-\eta^{i-1}}{h}\left[\lt(\frac{v^i-v^{i-1}}{\eta^i-\eta^{i-1}}\rt)' + \frac{(v^i-v^{i-1})^2/h^2}{(\eta^i-\eta^{i-1})^2/h^2}\right]\right|^2\\
&\le C\mathfrak{E}_N(0) + C\int_0^t \mathfrak{E}_N^2(s)ds.
\end{split}
\ee
\hfill $\Box$

For $h=1/N$, we define the functions $\eta^h(x,t)$ and $v^h(x,t)$, respectively, by
$$\eta^h(x,t)=\eta^{i,h}(x,t),\quad v^h(x,t)=v^{i,h}(x,t),\quad x\in (x_{i-1},x_i),\;1\le i\le N,\; t\in [0, T^*],$$
 where
\be\label{A29}
\begin{split}
\eta^{i,h}(x,t)=\eta^{i-1}(t)+ \frac{1}{h}(\eta^i(t)-\eta^{i-1}(t))(x-x_{i-1}),\\
v^{i,h}(x,t)=v^{i-1}(t)+ \frac{1}{h}(v^i(t)-v^{i-1}(t))(x-x_{i-1}).
\end{split}
\ee
 It follows that
\bee
\eta^h_t(x,t)=v^h(x,t),\quad \eta^h_x(x,t)= \frac{\eta^i(t)-\eta^{i-1}(t)}{h},\quad
 v^h_x(x,t)= \frac{v^i(t)-v^{i-1}(t)}{h},
\eee
for $x_{i-1}<x<x_i$, $1\le i\le N$ and $0\le t\le T^*$.
\begin{lem}\label{compactness} (Compactness)
  There exist functions $\eta(x,t)$ and $v(x,t)$, such that the subsequences of $\{ \eta^h(x,t)\}$ and $\{v^h(x,t) \}$, which is still denoted as $\{ \eta^h(x,t)\}$ and $\{v^h(x,t) \}$ for convenience, such that
  \bee
\eta^h\to \eta,\quad v^h\to v,\quad \eta^h_x \to \eta_x, \quad v^h_x\to v_x \quad a.e.\;  I\times[0,T^*].
  \eee
  Moreover, $\eta(x,t)$ and $v(x,t)$ satisfy $\eta_t=v$ and $\ef{6}$ almost everywhere, and the regularity in \ef{regu}.
\end{lem}
{\bf Proof.} From \ef{AA1}, the functions of the families $\{\eta^h\}$ and $\{v^h\}$ as functions of $x$ have uniformly bounded variations with respect to $h$ for each fixed $t\in[0,T^*]$. Suppose that the countable set $\{t_k\}_{k=1}^\infty$ is   dense in $[0,T^*]$. Then by Helly's theorem and a diagonal process, there exists a subsequence of $\{\eta^h\}$, which is still denoted as $\{\eta^h\}$ for simplicity, converging boundedly and almost everywhere for $x\in I$ and $\{t_k\}_{k=1}^\infty\subset [0,T^*]$ as $h\to 0.$ Employing Lebesgue's theorem, $\{\eta^h\}$ also converges in
$L^2(I)$ on $\{t_k\}_{k=1}^\infty$. Next, we show the continuity in time for the $L^2$-norm of $\{\eta^h\}$, that is,
\be\label{A30}
\begin{split}
\|\eta^{i,h}& (\cdot,t)-\eta^{i,h}(\cdot,s)\|_{L^2(I)}^2\\
 & \le  Ch\sum_{i=1}^{N-1} |\eta^i (t)-\eta^i(s)|^2+ C\sum_{i=1}^{N-1}\frac 1 h \int_s^t |v^i(\xi)-v^{i-1}(\xi)|^2 d\xi \\
 &\le  C\mathfrak{E}_N(0)|t-s|,
\end{split}
\ee
from the definition of $\eta^{i,h}$ in \ef{A29} and the uniform estimates \ef{AA1}. Therefore, due to the density of $\{t_k\}_{k=1}^\infty$ in $[0,T^*]$ and the continuity of $\|\eta^{i,h} (\cdot,t)\|_{L^2(I)}$ in time, there exists a function $\eta(x,t)$ such that $\{\eta^h\}$ converges to $\eta$ in $L^2(I)$ uniformly in $t\in [0,T^*]$. It follows that $\{\eta^h\}$ converges to $\eta$ $a.  e.$ in $x\in I$ and $t\in [0,T^*]$. Using similar arguments and the estimates \ef{AA1}, we find $\{v^h\}$, $\{v^h_x\}$ and $\{\eta^h_x\}$ also converges to $v$, $v_x$ and $\eta_x$, respectively, $a.  e.$ in $x\in I$ and $t\in [0,T^*]$. Since $\eta^h_t=v^h$, we get $\eta_t =v$ in $\mathcal{D}'(I\times[0,T^*])$ and  $a. e.$ in $ I\times [0,T^*]$. With the above convergence results and uniform estimates on the approximate sequences, we may verify that
 \be\label{A31}
\begin{split}
 \int_0^{T^*}\int_I &  \bar\rho v \phi_t  dxdt -\int_I \bar\rho v \phi  dx\big|_{t=0}^{T^*} = \int_0^{T^*}\int_I  \{\bar\rho^\gamma [1-(\eta_x)^{-\gamma}] + \eta_x^{-1} v_x\}   \phi_x  dxdt,
\end{split}
\ee
for any $\phi\in C^1([0,T^*],C^1_c (I))$. This shows that $\ef{6}_1$ is satisfied in the sense of distribution and almost everywhere.
Clearly, the regularity \ef{regu} of solutions  follow from the uniform bounds in Lemma \ref{uniform}, which can define the traces and boundary conditions of $(\eta,v)$. \hfill$\Box$

\begin{lem} There exists at most one solution $(\eta,v)$ to \ef{6} with $\eta_t=v$.
\end{lem}
{\bf Proof.} Suppose that there exist two solutions $(\eta_1,v_1)$ and $(\eta_2,v_2)$ to \ef{6} with $\eta_{it}=v_i$, $i=1,2$, satisfying the regularity \ef{regu} and the same initial and boundary conditions. Set $\zeta=\eta_1-\eta_2$ and $w=v_1-v_2$. Then $(\zeta,w)$ satisfies the following equations:
\be\label{B1}\lt\{\begin{split}
& \zeta_t=w \ \  & {\rm in} & \ \ I \times (0,T^*),\\
& \bar\rho w_t + \lt[  \bar\rho^\gamma (\eta_{1x}^{-\gamma}-\eta_{2x}^{-\gamma})\rt]_x = \lt(\frac{ {v_{1x}}}{\eta_{1x} }-\frac{ {v_{2x}}}{\eta_{2x} }\rt)_x  \ \  & {\rm in} & \ \ I \times (0,T^*),   \\
&w |_{x=0}=w_x  |_{x=1}=0,  \\
&(\zeta,w)|_{t=0}=(0,0).
\end{split}\rt.\ee
From $\ef{B1}_1$, we have, for $t\in[0,T]$,
\be\label{B4}
\int_I \zeta^2 dx (t)=2\int_0^t\int_I \zeta w dxdt\le \delta\int_0^t \int_I w^2 dxdt + C(\delta) \int_0^t \int_I \zeta^2 dxdt,
\ee
and similarly
\be\label{B5}
\int_I \zeta_x^2 dx (t) \le \delta\int_0^t \int_I w_x^2 dxdt + C(\delta) \int_0^t \int_I \zeta_x^2 dxdt.
\ee
Next, we multiply $\ef{B1}_2$ by $w$ and integrate over $I\times(0,t)$ to get
\be\label{B6}
\frac 1 2 \int_I \bar\rho  w^2 dx(t) \equiv J_1 +J_2,
\ee
where \bee\begin{split}
J_1&=- \int_0^t\int_I \lt(\frac{ {v_{1x}}}{\eta_{1x} }-\frac{ {v_{2x}}}{\eta_{2x} }\rt)w_x dx dt\\
&=- \int_0^t\int_I \frac{w_x^2}{\eta_{1x}}dxdt+\int_0^t\int_I \frac{v_{2x}}{\eta_{1x} \eta_{2x} }\zeta_x w_x dxdt\\
   &\le - \frac 1 2 \int_0^t\int_I \frac{w_x^2}{\eta_{1x}}dxdt +C  \int_0^t\int_I \frac{v_{2x}^2}{\eta_{1x}\eta_{2x}^2}\zeta_x^2 dxdt,
\end{split}\eee
and
\bee\begin{split}
J_2&=  \int_0^t\int_I \bar\rho^{\gamma}\lt( \eta_{1x}^{-\gamma}-\eta_{2x}^{-\gamma}\rt)w_x dx dt\\
&\le \delta \int_0^t\int_I  \frac{w_x^2}{\eta_{1x}}  dxdt+  C(\delta)\int_0^t\int_I \eta_{1x}(\eta_{1x}^{-2\gamma}+\eta_{2x}^{-2\gamma})\zeta_x^2    dxdt.
\end{split}\eee
It follows from the fact that $\eta_{ix},\eta_{ix}^{-1}, v_{2x}\in L^\infty_{x,t}$, $i=1,2$,
\be\label{B7}
\frac 1 2 \int_I \bar\rho  w^2 dx(t) + \frac 1 4 \int_0^t\int_I  \frac{w_x^2}{\eta_{1x}}  dxdt\le C\int_0^t\int_I \zeta_x^2    dxdt,\q t\in[0,T^*].
\ee
Combining \ef{B4}, \ef{B5} and \ef{B7}, using  Gronwall's inequality, Lemma \ref{lem21} and the fact that $\eta_{1x}^{-1}$ is bounded from below, we show this lemma by choosing $\delta$ to be small. \hfill$\Box$

\section*{Acknowledgements.} This work was partially supported by NSFC under grant 11971477 and 11761141008, the Fundamental Research Funds for the Central Universities and the Research Funds of Renmin University of China (grant No. 18XNLG30). This work was initiated when the author visited the City University of Hong Kong. The author would like to thank Prof. Xianpeng Hu for helpful discussions and hospitality. The final version of this paper was completed when the author visited Brown University under the support of the China Scholarship Council (grant No. 201806365010). Moreover, I would like to thank the anonymous reviewers for their valuable suggestions and comments, which significantly help me to improve the presentation of this paper.

\end{document}